\documentclass[12pt]{amsart}

\usepackage{amsmath,amsthm,amsfonts,latexsym,amssymb,amscd,enumerate,times,fullpage}
\usepackage[all]{xy}
\newtheorem{theorem}{Theorem}[section]
\newtheorem{thm}[theorem]{Theorem}
\newtheorem{prop}[theorem]{Proposition}
\newtheorem{lem}[theorem]{Lemma}
\newtheorem{cor}[theorem]{Corollary}
\newtheorem{conj}[theorem]{Conjecture}

\theoremstyle{remark}

\theoremstyle{definition}
\newtheorem{defn}[theorem]{Definition}

\theoremstyle{definition}

\theoremstyle{remark}

\DeclareMathOperator{\im}{im}      

\DeclareMathOperator{\const}{{\rm const}}

\DeclareMathOperator{\cone}{{\rm cone}}

\DeclareMathOperator{\Mor}{Mor} 

\newcommand{\botimes}{\, \overline \otimes \,}

\newcommand{\coker}{{\rm coker}\,} 
\newcommand{\F}{\mathbb{ F}}
\newcommand{\N}{\mathbb{ N}}

\newcommand{\Z}{\mathbb{ Z}}

\newcommand{\Hom}{\operatorname{Hom}}

\newcommand{\Q}{\mathbb{ Q}}
\newcommand{\R}{\mathbb{ R}}

\newcommand{\id}{\operatorname{id}}

{\catcode`@=11\global\let\c@equation=\c@theorem}

\date{\today; MSC 2000: primary 57S17; secondary 55P62}
\begin{document}

\title{The stable free rank of symmetry of products of spheres} 

\author{Bernhard Hanke}

\begin{abstract} A well known conjecture in the theory of transformation groups
states that if $p$ is a prime and $(\Z/p)^r$ 
acts freely on a product of $k$ spheres, then $r \leq k$. We prove 
this assertion if $p$ is large compared 
to the dimension of the product of spheres. The argument builds 
on tame homotopy theory for non-simply connected spaces. 
\end{abstract} 

\address{Universit{\"a}t M{\"u}nchen \\ Germany \\}
\email{hanke@mathematik.uni-muenchen.de \vspace{1cm}}

\maketitle

\section{Introduction} 
Transformation group theory investigates 
symmetries of topological spaces. An important aspect of 
this program is to define and study invariants that distinguish spaces 
admitting 
lots of symmetries from less symmetric ones. 

In this paper we concentrate on 
one of these invariants, the so called {\em free $p$-rank}
\[
       \max \{ r \, | \,  (\Z/p)^r {\rm~acts~freely~on~} X \} \, ,  
\]
defined for any topological space $X$ and any prime number $p$ . Here all groups are acting topologically. 
We recall the following fact from classical Smith theory. 

\begin{thm} \label{smith} The free $p$-rank of $S^n$ is equal to $\begin{cases} 1  {\rm~for~odd~}n  \\
                            1  {\rm~for~even~}n{\rm~and~} p=2  \\
                            0  {\rm~for~even~}n{\rm~and~} p>2 \, . 
              \end{cases} $
\end{thm} 

In view of this theorem it is natural to look for a corresponding 
result, if $X$ is not just a single sphere, but a product of spheres,  
\[
    X = S^{n_1} \times S^{n_2} \times \ldots \times S^{n_k} \, . 
\]
The following statement appears in several places in the literature either 
as a question \cite[Question 7.2]{AB}, \cite[Problem 809]{MR} or as a conjecture 
\cite[Conjecture 2.1]{A}, \cite[Conjecture 3.1.4]{AD}. 

\begin{conj} \label{verm} If $(\Z/p)^r$ acts freely on $X$, then $r \leq k$. 
\end{conj}

Actually, if $p$ is odd, it is (in view of Theorem \ref{smith}) reasonable to conjecture 
that $r$ is bounded above by the number $k_o$ of odd dimensional spheres in $X$. Conjecture \ref{verm}, 
in this sharper form for odd $p$, has been verified in the following cases:  

\begin{enumerate}[1)]
   \item $k \leq 2$, see Heller \cite{He};  $k\leq 3$, $p=2$,  see Carlsson \cite{Ca6}. 
   \item $n_1 = \ldots = n_k$ and  in addition  
       \begin{enumerate}[a)]
           \item the induced action on integral homology is trivial, see Carlsson \cite{Ca1}, or
           \item the induced action on integral homology is unrestricted, but 
                 if $p=2$, then $n_i \neq 3,7$, see Adem-Browder \cite{AB} (for $p \neq 2$ or $n_i \neq 1,3,7$),
                 and Yal\c c\i n \cite{Y} (for $p=2$ and $n_i=1$).
       \end{enumerate}
   \item assume $n_1 \leq n_2 \leq \ldots \leq n_k$. Then $n_1 \geq 2$ and 
         for all $1 \leq i \leq k-1$ either $n_i = n_{i+1}$ or 
         $2n_i \leq n_{i+1}$. Furthermore, the induced $(\Z/p)^r$-action 
         on $\pi_*(X)$ is trivial and $p > 3 \dim X$.  See S\"orensen \cite{So}. 
  \end{enumerate}

The following theorem is our main result. It settles a stable form of Conjecture \ref{verm}. 

\begin{thm} \label{main} If $p > 3 \dim X$  and $(\Z/p)^r$ acts freely on $X$, then $r \leq  k_o$\, . 
\end{thm}

Because $(\Z/p)^{k_o}$ acts freely on a product of $k_o$ odd dimensional spheres, this result implies that the free $p$-rank 
of $X$ is equal to $k_o$ assumed that $p > 3 \dim X$. Theorem \ref{main} also leads  
to the following estimate of the {\em free toral rank} of $X$ which follows
from \cite[Theorem T]{Hal}. 

\begin{thm} \label{rational} If $(S^1)^r$ acts freely on $X$, then $r \leq k_o$\, .
\end{thm}

However, even for very large primes Theorem \ref{main} cannot be deduced from this 
result by some sort of limiting process: Browder \cite{B2} constructed free 
actions of $(\Z/p)^r$ on $(S^m)^k$ for each odd $m \geq 3$ and 
$r \geq 4$, $k \geq r$, $p>km/2$ which are exotic in the sense that they cannot be 
extended to $(S^1)^r$-actions. 

Theorem \ref{rational} is implied by a theorem of Halperin \cite{Hal} which yields upper bounds of 
free toral ranks in terms of homotopy Euler characteristics. The
proof is based on rational homotopy theory applied to the Borel space $X_{(S^1)^r}$.  
Hence rational homotopy theory has become an important tool in the study of 
torus actions in general, see \cite[Section 2]{AP} for a survey. 

Because the rational homotopy type of $B(\Z/p)^r$ is trivial, it is obvious that in the context of 
Theorem \ref{main}, rational homotopy theory cannot be applied in a reasonable way to study 
the Borel space $X_{(\Z/p)^r}$. But {\em tame homotopy} 
theory seems a more promising approach. This theory was invented 
by Dwyer \cite{Dw} and is modelled on  Quillen's rational homotopy theory \cite{Qu}, 
but without immediately losing $p$-torsion information for all primes $p$. 
One of the first observations in tame homotopy theory
may be phrased as follows: if $X$ is an $(r-1)$-connected space ($r \geq 1$)
and $\pi_{r+k}(X)$ is a $\Z[p^{-1} \,  | \, 2p-3 \leq k]$-module for every $k \geq 1$, then 
the complexity of the Postnikov invariants of $X$ should be comparable to that of a rational 
space due to the vanishing of the relevant higher reduced Steenrod power operations in the Postnikov pieces of $X$.

Nevertheless, in the original setup, Dwyer's theory could only be formulated for $r > 2$, 
see \cite[1.5]{Dw}. In the tame setting, the restriction to 
simply connected spaces, which excludes the Borel spaces for $(\Z/p)^r$-actions, 
seems unavoidable due to the non-vanishing of the 
$p$-fold Massey products $H^1(B \Z/p;\F_p) \to H^2(B \Z/p ; \F_p)$ for any prime $p$, see \cite[Theorem 14]{Kraines}.
We will make more comments on this point later on.

To any space $X$  one can associate the {\em Sullivan-de Rham algebra} \cite{BG, Sul}, 
a commutative graded differential algebra over $\Q$ modelled on de Rham differential forms, 
which calculates the rational cohomology of $X$. This construction widens the scope of rational 
homotopy theory from simply connected to nilpotent spaces. A distinctive feature 
of this approach is the construction of a {\em minimal model} for any space $X$ out of 
the Sullivan-de Rham algebra. The minimal model still calculates the rational cohomology of the 
underlying space, but in addition its associated vector space of indecomposables 
can be identified with the dual of $\pi_{*}(X) \otimes \Q$, if  $X$ is nilpotent, $H^*(X;\Q)$ is finite dimensional 
in each degree and  $\pi_1(X)$ is abelian. 

For a proof of Theorem \ref{rational} one writes down a minimal 
model of the Borel space $X_{(S^1)^r}$ and argues that if $r > k_o$, then the cohomology of 
the minimal model is nonzero in arbitrarily high degrees. But this contradicts the fact that by the freeness of the 
action, $X_{(S^1)^r}$  is homotopy equivalent to the finite dimensional space $X/(S^1)^r$. 
The analysis of the minimal model that leads to these conclusions is carried out in the 
fundamental paper \cite{Hal}.

The generalization of the Sullivan-de Rham algebra to the tame 
setting is realized in the work of Cenkl and Porter \cite{CP} and 
is achieved by considering commutative graded differential algebras over $\Q$ which are 
equipped with {\em filtrations} as a new structural ingredient. By definition, elements in filtration 
level $q$ are divisible by any prime $p \leq q$, but not necessarily by larger primes, and 
filtration levels  are additive under multiplication of elements. The construction of 
the Cenkl-Porter complex is based on differential forms similar to the Sullivan-de Rham algebra, but  
in order to keep sufficient control over filtration levels one  works with differential forms defined on a cubical decomposition 
of the standard simplex. 
The {\em Cenkl-Porter theorem} states that the integration of forms yields 
a cochain map of the Cenkl-Porter complex in filtration level $q$ to the singular 
cochain complex  with $\Z[p^{-1} \, | \, p \leq q ]$-coefficients which in cohomology induces multiplicative isomorphisms. 
We will recall the essential steps of this construction 
in Section \ref{CenklPorter} of our paper.  Our proof of the Cenkl-Porter theorem is different from the original approach 
in \cite{CP} in that we do not analyse the integration map, which is not multiplicative on the cochain level, but 
construct multiplicative cochain maps inducing isomorphisms in cohomology from the Cenkl-Porter and 
singular cochain complexes to a third cochain complex. This approach is inspired by the discussion 
of the de Rham theorem in \cite[Theorem II.10.9]{FHT} and enables us to take up and resolve 
the issue of  possibly non-vanishing higher Massey products at the end of Section \ref{CenklPorter}. 

It is now desirable to replace the Cenkl-Porter complex by a smaller filtered commutative graded 
differential algebra which nicely reflects the homotopy type of the underlying space 
in a similar manner as the minimal algebra does in rational 
homotopy theory. In the tame setting, the most transparent and elementary
approach to the construction of such an algebra is based on  a {\em tame Hirsch lemma} which 
is used to build the desired small filtered cochain 
algebra along a Postnikov decomposition of the underlying space, compare the discussion in \cite{DGMS} for the rational case. 
For simply connected spaces 
such a result is proven in \cite[p. 203]{Sch}. The non-simply connected case, 
which is relevant for our purposes,  is much more involved and was carried out in the remarkable 
diploma thesis of S\"orensen \cite{So}. After introducing the necessary notions we will  
state the tame Hirsch lemma in Section \ref{Hirsch}. The proof of this 
result, which is not available in the published literature,  will be given  
in an appendix to our paper along the lines of \cite{So}. I owe my thanks 
to Till S\"orensen, who left academia some seventeen years ago, for producing a PDF scan 
of his work and to the FU Berlin for storing it on their preprint server.

With this machinery in hand we construct small approximative commutative $\F_p$-cochain models 
of Borel spaces associated to $(\Z/p)^r$-spaces $X$ under fairly general assumptions  
in Section \ref{small}. The main 
idea is to use a simultaneous Postnikov decomposition of the fibre and total space of the 
Borel fibration $X \hookrightarrow X_{(\Z/p)^r} \to B(\Z/p)^r$.

The cochain model resulting from this discussion is the starting point for our proof of 
Theorem \ref{main} in Section \ref{cool}. The argument is inspired in part by the 
paper \cite{Hal}, which provides techniques to simplify the analysis of free commutative
graded differential algebras over fields of characteristic zero. 
We emphasize that as stated, some of 
the arguments in {\em loc.~cit.} do not work over fields of prime 
characteristics. However, in the context of actions on 
products of spheres, with which we are mainly concerned, the    
resulting $\F_p$-cochain algebras are special
enough so that the proof of Theorem \ref{main} can be completed fairly quickly by  a direct 
argument.   

The main result of Section \ref{small}, Theorem \ref{smallcochain}, may be useful 
for other purposes in the cohomology theory of $p$-torus actions, which is 
now closely tied to the cohomology theory of $(S^1)^r$-actions, if $p$ is large. 
For example we hope that it is possible to generalize \cite[Theorem T]{Hal} in 
the sense that for large $p$, the free $p$-rank of a finite connected CW-complex with abelian fundamental 
group is bounded above by the rational homotopy Euler characteristic of $X$. 
However, we will not discuss this in greater detail.

We believe that the methods developed in this paper are not sufficient
to establish the general form of Conjecture \ref{verm} for small primes. This problem
remains open. 

Acknowledgements. I am grateful to Bill Browder for drawing  my attention to \cite{So}, 
to Jonathan Bowden for proofreading a previous version of this manuscript and to 
the referees for numerous helpful comments and remarks.

\section{Proof of Theorem \ref{main}} \label{cool} 

Let $p$ be a prime and assume $G := (\Z/p)^r$ acts freely on 
$X= S^{n_1} \times \ldots \times S^{n_k}$. 
Denoting by $k_e$ the number of even dimensional spheres in $X$, we assume that the 
dimensions $n_1, \ldots, n_{k_e}$ are even 
(and greater than $0$) and $n_{k_e+1}, \ldots, n_{k}$ are odd. 
Because $G$ acts freely, the Borel space  $X_G := EG \times_G X$
is homotopy equivalent to the orbit space $X/G$. This is a topological manifold of 
dimension $\dim X$ and consequently $H^n(X_G ; \F_p) = 0$ 
for all $n > \dim X$.  Assume that $p > 3 \dim X$. In particular, $p$ is odd.

Under these assumptions the following result will be proved at the end of Section \ref{small}.

\begin{prop} \label{later} 
There are free associative, commutative graded differential algebras with unit  
$(E^*,d_E)$ and $(M^*,d_M)$ over $\F_p$ with the following 
properties:  \begin{enumerate}[1)]
\item as graded algebras $E^* =   \F_p[t_1, \ldots, t_r] \otimes \Lambda(s_1, \ldots, s_r) \otimes M^*$, where 
         each $t_i$ is of degree $2$ and each $s_i$ is of degree $1$. 
\item the differential $d_E$ is zero on $\F_p[t_1, \ldots, t_r] \otimes \Lambda(s_1, \ldots, s_r) \otimes 1$ and the map 
         $E^* \to M^*$, $t_i, s_i \mapsto 0$, is a cochain map. 
\item $M^*$ is the tensor product of commutative graded differential algebras $(M^*_j, d_{M_j})$, $1 \leq j \leq k$, that
         are defined as follows: 
  \begin{enumerate}[a)] 
   \item  $M_j^*=\Lambda(\sigma_j)$ with $\deg(\sigma_j) = n_j$ and $d_{M_j}(\sigma_j) = 0$, if  $j > k_e$ (i.e.~if $n_j$ is odd), 
   \item  $M_j^*=\F_p[\tau_j]$ with $\deg(\tau_j)   = n_j$ and  $d_{M_j}(\tau_j) = 0$, if $j \leq k_e$ and $2n_j-1 > \dim X +1$\, , 
   \item  $M_j^*=\F_p[\tau_j] \otimes \Lambda(\eta_j)$ with  $\deg(\tau_j) = n_j, \deg(\eta_j) = 2n_j -1$ and 
             $d_M(\tau_j) = 0, d_M(\eta_j) = \tau_j^2$, if  $j \leq k_e$ and $2n_j -1 \leq \dim X +1$, 
\end{enumerate}
\item each monomial in $t_1, \ldots, t_r$ of cohomological degree at least $\dim X +1 $ represents the zero class in $H^*(E)$. 
\end{enumerate}
\end{prop}

The short exact sequence of $\F_p$-cochain algebras 
\[
    \F_p[t_1, \ldots, t_r] \otimes \Lambda(s_1, \ldots, s_r) \to E^* \to M^*
\]
models the Borel fibration $X \hookrightarrow X_G \to BG$ up to degree $\dim X + 1$ (respectively up to degree $\dim X + 2$ 
as far 
as the monomials in $t_1, \ldots , t_r$ are concerned) implying the last statement of Proposition \ref{later}. 
In rational homotopy theory the Grivel-Halperin-Thomas
theorem \cite[Theorem 2.5.1]{AP} yields a similar short exact sequence of rational differential graded algebras 
for spaces equipped with 
$(S^1)^r$-actions that models the Borel fibration in all degrees. In this case $\F_p[t_1, \ldots, t_r] 
\otimes 
\Lambda(s_1, \ldots, s_r)$ is replaced by a polynomial algebra $\Q[t_1, \ldots, t_r] $, the minimal model 
of $B(S^1)^r$, and  the algebra $M^*$ can be 
chosen as the minimal model of $X$. If $X$ is a product of spheres as before, this minimal model looks very 
similar to the algebra $M^*$ 
appearing in Proposition \ref{later}. In Section \ref{small} tame homotopy theory will be used to establish a 
counterpart to the Grivel-Halperin-Thomas theorem 
for $G$-spaces and cohomology with $\F_p$-coefficients which models the Borel fibration up to degree $\dim X + 1$ 
and leads  to the statement of Proposition \ref{later}. 
In  the description of $M^*$ the two cases $2 n_j -1 > \dim X + 1$ and $2 n_j -1 \leq \dim X +1$ need to be distinguished, 
because for $p > 3 \dim X$ the procedure 
of Section \ref{small} constructs free generators  of 
the algebra $E^*$ only up to degree $\dim X + 1$. When the algebra $E^*$ has been simplified, 
we can easily add generators $\eta_j$  of degrees $2n_j - 1 > \dim X +1$, see below.

Starting from $(E^*,d_E)$ we shall construct a commutative graded differential algebra $(F^*,\delta_F)$ that 
is easier to handle.  
Let 
\[
    F^* := \F_p[t_1, \ldots, t_r] \otimes M^* 
\]
be obtained from $E^*$ by dividing out 
the ideal $(s_1, \ldots, s_r)$. The induced differential on $F^*$ is denoted 
by $d_F$. 

Inspired by the construction of {\em pure towers} in \cite{Hal}, 
we deform $d_F$ to  another differential $\delta_F$ on $F^*$ as follows: 
$\delta_F$ is a derivation that vanishes on $\F_p[t_1, \ldots, t_r, \tau_1, \ldots, \tau_{k_e}]$ and satisfies 
\begin{eqnarray*} 
  \delta_F(\sigma_j) & = & \pi(d_F(\sigma_j))   \\
  \delta_F(\eta_j)   &  =  &  \pi ( d_F(\eta_j)) \, .   
\end{eqnarray*} 
where $\pi: F^* \to F^*$ is the projection onto 
$\F_p[t_1, \ldots, t_r, \tau_1, \ldots, \tau_{k_e}]$ given by evaluating the 
odd degree generators $\sigma_j, \eta_j$ at $0$. 
It is easy to verify that $\delta_F^2 = 0$. 

For $l \geq 0$ let $\Sigma^l \subset F^*$ be the $\F_p[t_1, \ldots, t_r]$-linear
subspace generated by the monomials in $M^*$ containing exactly $l$ of the  odd degree generators $\sigma_j, \eta_j$, 
in particular $\Sigma^{l} = 0$ for $l > k$ by the anticommutativity of 
the product (recall that $p$ is odd). We set $\Sigma^+ := \bigoplus_{l \geq 1} \Sigma^l$. 
This is a nilpotent ideal in $F^*$. 

\begin{lem} \label{easy} For all $l \geq 1$ the differential $\delta_F$ maps $\Sigma^l$ to $\Sigma^{l-1}$. Furthermore, 
the image of  $\delta_F - d_F$ is contained in $\Sigma^+$. 
\end{lem} 

\begin{proof} The first assertion holds by the definition and derivation property of $\delta_F$. 

The second assertion holds for the generators $\sigma_j$ and $\eta_j$, because $\im (\id - \pi) \subset \Sigma^+$, 
it holds for the generators $t_i$, because $\delta_F$ and $d_F$ send these elements 
to zero and it holds for the generators $\tau_j$, because each $d_F(\tau_j)$ is  of odd degree and therefore contained in 
$\Sigma^+$. 
This implies the second assertion in general, since $\Sigma^+$ is an ideal in $F^*$ and $\delta_F - d_F$ is a derivation. 
\end{proof}

The elements $t_i$, $1 \leq i \leq r$,  and $\tau_j$, $1 \leq j \leq k_e$, represent cocycles 
in $(F^*, \delta_F)$. Let $[t_i]$ and $[\tau_j]$ be the corresponding cohomology classes. 
The proof of Theorem \ref{main} depends on the following crucial fact. 

\begin{prop} \label{genial} The elements $[t_i]$ are nilpotent in $H^*(F,\delta_F)$. The same 
holds for $[\tau_j]$, if $2n_j-1 \leq \dim X +1$. 
\end{prop} 

\begin{proof} As noted earlier 
each monomial in $t_1, \ldots, t_r$ of cohomological degree at least $\dim X +1 $ represents the zero class in $H^*(E)$. 
Let $m$ be 
such a monomial and write $m = d_E (c)$ for a cochain $c \in E^*$. We will show that $m$ is a 
coboundary in $(F^*,\delta_F)$ as well. 

First, after dividing out the ideal $(s_1, \ldots, s_r)$, we get an 
analogous  equation $m = d_F (c)$ for some $c \in F^*$. By Lemma \ref{easy} 
we conclude  $\delta_F (c) = m  + \omega$
where $\omega \in \Sigma^+$. Let $c_1$ be the component of $c$ in $\Sigma^1$.  
Lemma \ref{easy} and the fact that $m \in \Sigma^0$ imply the equation $\delta_F(c_1) = m$.  
This shows that $m$ is a coboundary in $(F^*,\delta_F)$. 

In particular we have shown  
that the classes $[t_i] \in H^*(F,\delta_F)$ are nilpotent. 

The cochain algebra $(F^*,\delta_F)$ has a decreasing filtration given by 
\[
   \mathcal{F}_{\gamma}(F^*) =  \F_p[t_1, \ldots, t_r]^{\geq \gamma} \otimes M^*  
\]
where $\gamma \in \N$ denotes the cohomological degree. Our previous argument and the fact that each $\tau_j$ is 
a cocycle in $(F^*, \delta_F)$ imply that each element in  $\Sigma^0 \subset F^*$ 
in filtration level at least $\dim X + 1$ is a coboundary in $(F^*,\delta_F)$.  

Now pick a  $j \in \{1, \ldots, k_e\}$ with $2n_j -1 \leq \dim X +1$. Recall (see the 
description of $E^*$ in Proposition \ref{later}) that
\[  
    d_F ( \eta_j ) =  \tau_j^2 \mod \mathcal{F}_{2}(F) \, . 
\]
By the definition of $\delta_F$ we obtain 
\[ 
    \delta_F(\eta_j) = \pi(\tau_j^2) = \tau_j^2 \mod   \mathcal{F}_2(F) \, , 
\]
since the map $\pi$ preserves the ideal $(t_1, \ldots, t_r)=\mathcal{F}_2(F)$. 
This implies that $\tau_j^2$ 
is $\delta_F$-cohomologous to a cocycle $c \in \mathcal{F}_{2}(F^*,\delta_F)$. 
Hence, $(\tau_j^{2})^{\dim X}$ is $\delta_F$-cohomologous to $c^{\dim X} \in \mathcal{F}_{2 \dim X}(F)$.
We can split $c^{\dim X}$ into a sum $c_0 + c^+$ where  $c_0 \in \Sigma^0 \cap \mathcal{F}_{2 \dim X}(F)$ 
and $c^+ \in \Sigma^+ \cap \mathcal{F}_{2 \dim X}(F)$. 
As noted earlier, $c_0$ is $\delta_F$-cohomologous to zero. 
Because $\Sigma^+$ is nilpotent, the element $c^+$ is 
nilpotent. 

We conclude that $\tau_j^{2 \dim X}$ is $\delta_F$-cohomologous to a nilpotent cocycle in $(F^*,\delta _F)$. 
\end{proof} 

After these preparations we can finish the proof of Theorem \ref{main}. 
At first we adjoin new exterior generators $\eta_j$ of degree $2n_j -1$ to  $F^*$ for 
all $j \in \{1, \ldots, k_e\}$ satisfying $2n_j -1 > \dim X +1$. The differential $\delta_F$ 
is extended to a differential on this new graded commutative algebra by setting 
$\delta_F(\eta_j) := \tau_j^2$, if $2n_j -1 > \dim X +1$. This is 
possible, because $\delta_F(\tau_j) = 0$ by construction of $\delta_F$. 
This new commutative graded differential algebra is still denoted by $(F^*, \delta_F)$.
Together with  Proposition \ref{genial} we see that 
the elements $t_i$, $1 \leq i \leq r$, and $\tau_j$, $1 \leq j \leq k_e$, 
define nilpotent classes in $H^*(F,\delta_F)$. This implies that  $H^*(F, \delta_F)$ is a finite dimensional 
$\F_p$-vector space. 

We consider the ideal  
\[
  I = \big( \delta_F (\eta_1) , \ldots, \delta_F(\eta_{k_e}), \delta_F(\sigma_{k_e+1}) , \ldots, \delta_F(\sigma_{k}) \big) \subset 
       \F_p[t_1, \ldots, t_r, \tau_1, \ldots, \tau_{k_e} ] 
\]
contained in $\im(\delta_F)$ and obtain an inclusion 
\[
   \F_p[t_1, \ldots, t_r, \tau_1, \ldots, \tau_{k_e}] / I \subset H^*(F,\delta_F) \, . 
\]
Here we use the fact that the coboundaries in $(F^*,\delta_F)$ are contained in the ideal 
$I \cdot F^*$, whose intersection with $\F_p[t_1, \ldots, t_r, \tau_1, \ldots, \tau_{k_e}]$ is equal to $I$.  This 
can be seen by applying the evaluation map $\pi$. 
We conclude that $\F_p[t_1, \ldots, t_r, \tau_1, \ldots, \tau_{k_e}] / I$ is a finite dimensional 
$\F_p$-vector space. Because $I$ is generated by homogenous elements of positive degree, it does not contain a unit 
of $\F_p[t_1, \ldots, t_r,\tau_1, \ldots, \tau_{k_e}]$ and hence there is a minimal prime 
ideal  $\mathfrak{p} \subset \F_p[t_1, \ldots, t_r, \tau_1, \ldots, \tau_{k_e}]$ containing $I$. The quotient 
$\F_p[t_1, \ldots, t_r, \tau_1, \ldots, \tau_{k_e}] / \mathfrak{p}$ is both a finite dimensional $\F_p$-vector
space and an integral domain. Hence $\mathfrak{p} = (t_1, \ldots , t_r, \tau_1, \ldots, \tau_{k_e})$ and 
consequently ${\rm height} (\mathfrak{p}) = r + k_e$. By Krull's Principal Ideal Theorem, see \cite[Theorem 10.2]{Eisen}, 
the number of generators of $I$ must be at least  $r+k_e$. 
From the definition of $I$ we derive the inequality $k_e + k_o \geq r+k_e$. This 
implies $k_o \geq r$ and finishes the proof of Theorem \ref{main}.

\section{Tame homotopy theory via differential forms} \label{CenklPorter} 

The remainder of our paper is devoted to the construction of small approximative $\F_p$-cochain models of Borel 
spaces associated to $(\Z/p)^r$-spaces in Theorem \ref{smallcochain}. Then the corresponding cochain 
models for $(\Z/p)^r$-actions on products of spheres will lead directly to 
the cochain algebras $(E^*,d_E)$ and 
$(M^*,d_M)$ appearing in Proposition \ref{later}. 

First, in this section, we collect some fundamental notions and constructions 
from tame homotopy theory via polynomial forms as initiated in \cite{CP} and 
further developed in \cite{Bo,Ma,Ou,Sch,So,St}. 
As in the rational 
situation \cite{BG} it is convenient to formulate the theory in the language 
of simplicial sets. The reader unfamiliar with this may consult the introductory text \cite{May}.

Prime numbers are henceforth denoted  by 
the letter $l$.
For $q \in \N$ we set 
\[
    \Q_q := \Z [ l^{-1} \, | \, l \leq q ]\, , 
\]
i.e. $\Q_q$ is the smallest subring of $\Q$ where all primes $l \leq q$ are 
invertible. In particular, $\Q_0 = \Q_1 = \Z$.

\begin{defn} \label{basik} A {\em filtered cochain complex}  is a cochain complex $(V^*, d)$ over $\Q$, graded over $\N$, 
together with an increasing sequence of subcomplexes
\[
    V^{*,0} \subset V^{*,1} \subset V^{*,2} \subset \ldots \subset V^{*} 
\]
so that $V^{*,q}$ is a cochain complex of $\Q_q$-modules.  We call $V^{*,q}$ 
the subcomplex of {\em filtration level $q$}. 
A {\em filtered graded differential algebra} is a 
filtered cochain complex $(A^*,d)$ which is equipped with an associative 
multiplication with unit $1 \in A^{0,0}$, which restricts  to 
maps $A^{*, q} \otimes  A^{*,q'} \to A^{*, q+q'}$. Furthermore, $d$ is assumed to act as a graded derivation.  If the
multiplication is graded commutative, then $A^*$ is called a {\em filtered commutative graded differential algebra} 
or {\em filtered CGDA}. 
\end{defn}

If $A^{*}$ is a filtered cochain complex or filtered graded differential algebra and $x \in A^{p}$, 
we call $p$ the {\em degree} of  $x$ 
and $\min \{q \, | \, x \in A^{p,q} \}$ the {\em filtration degree} of $x$. The category 
of filtered CGDAs and filtration preserving CGDA maps is denoted by $\mathcal{A}$. 
It can be equipped with the structure 
of a closed model category \cite{Ma,Qu,Sch,So,St}, but this fact will  not be needed for 
our discussion. 

Let $V_1$ and $V_2$ be two filtered cochain complexes. Their {\em tensor product} is the 
cochain complex $V_1 \otimes V_2$ with the usual grading and differential together 
with the filtration  
\[
   (V_1 \otimes V_2)^{q} :=  \sum_{q_1 + q_2 = q} 
                   V_1^{*,q_1} \otimes V_2^{*,q_2} \otimes \Q_q \subset V_1^* \otimes V_2^*   \, , 
\]
where $\sum$ denotes a sum which need not be direct. The same definition applies to filtered CGDAs.   The tensor 
product defines the coproduct in $\mathcal{A}$.

If $V$ is a filtered cochain complex, let $\Lambda V$ denote the free CGDA generated by $V$. It is 
equipped with a filtration as described in the preceding paragraph. Note 
that this is an algebra with unit $1 \in (\Lambda V)^{0,0}$ and that  the filtration degree of 
a monomial $v_1 \cdot \ldots \cdot v_n$, $v_i \in V$,  is less than or equal to the sum of 
the filtration degrees of the $v_i$. 

Sometimes the  following variaton of the tensor product will be useful.  The {\em filtrationwise tensor product} $V_1 \botimes V_2$
is defined as the tensor product complex $V_1 \otimes V_2$, but equipped with the filtration 
\[
    (V_1 \botimes V_2)^{*,q} := V_1^{*,q} \otimes V_2^{*,q} \, . 
\]

The prototypical example of a filtered CGDA is provided by the cubical cochain algebra
of Cenkl and Porter \cite{CP} which refines the Sullivan-de Rham algebra \cite[\S 1]{BG} used in rational 
homotopy theory. Whereas the latter is based on forms on the standard
simplices $\Delta^{n} \subset \R^{n+1}$, the Cenkl-Porter construction works with 
forms on a cubical decomposition of $\Delta^n$. This is essential for 
defining a filtration on the resulting cochain algebra in such a way 
that in each positive filtration level one gets a cohomology theory in the sense of Cartan \cite{Ca}.  
We recall the essential steps of this construction. 

We cubically decompose the standard $n$-simplex by regarding it as the subset   
\[
    \Delta^n := \{ (x_0, \ldots, x_n) \in \R^{n+1} ~|~ 0 \leq x_i \leq 1, \prod_{i=0}^n x_i = 0 \}  
\]
of the boundary of the $(n+1)$-dimensional cube.  
The vertices $v_0, \ldots, v_n$ of $\Delta^n$ are given by
\[
    v_i := (1,1,\ldots, 1,0,1, \ldots, 1) 
\]
with $0$ located at the $i$-th entry, $0 \leq i \leq n$. The inclusion of the $i$-th face into $\Delta^n$
and the $i$-th collapse onto $\Delta^n$ ($0 \leq i \leq n$) are the maps 
\begin{eqnarray*}
  \delta_i : \Delta^{n-1} \to \Delta^n  &,  & (x_0, \ldots, x_{n-1}) \mapsto (x_0, \ldots, x_{i-1}, 1, x_{i+1}, \ldots, 
                                                                              x_{n-1})   \\ 
    \sigma_i : \Delta^{n+1} \to \Delta^{n} & ,  &  (x_0 , \ldots, x_{n+1}) \mapsto (x_0, \ldots, x_{i-1}, x_i \cdot x_{i+1}, 
                                                                             x_{i+2}, \ldots, x_{n+1}) \, . 
\end{eqnarray*}
We consider the free commutative graded algebra over $\Z$ 
\[
    \Z[t_0, \ldots, t_n] \otimes \Lambda_{\Z}(dt_0, \ldots, dt_n) 
\]
with generators $t_0, \ldots, t_n$ in degree $0$ and $dt_1, \ldots, dt_n$ in degree $1$. 
To a monomial 
\[
       t_0^{\alpha_0} dt_0^{\epsilon_0} \cdot \ldots \cdot t_n^{\alpha_n} dt_n^{\epsilon_n} \, , ~~ \alpha_i \geq 0 \, , ~~ 0 \leq \epsilon_i \leq 1 \, ,  
\]
in this algebra we assign the  
{\em filtration degree} $\max_i \{ \alpha_i + \epsilon_i\}$. 
The monomials of filtration degree $0$ are exactly the constant ones.
Let $I$ be the homogenous ideal generated by the monomials 
with $\alpha_i + \epsilon_i > 0$ for all $i$. 
The quotient 
\[
         \Z[t_0, \ldots, t_n] \otimes \Lambda_{\Z}(dt_0, \ldots, dt_n) / I
\]
is again a graded algebra and is called  the {\em algebra of compatible polynomial forms on the cubical 
decomposition of $\Delta^n$}. 

For $p,q \geq 0$, we define 
\[
     T_n^{p,q}(\Z)\subset (\Z[t_0, \ldots, t_n] \otimes \Lambda_{\Z}(dt_0, \ldots, dt_n) / I)^p
\]
as the $\Z$-module generated by the 
(cosets of) monomials of degree $p$ and filtration degree at most $q$. Exterior product of 
polynomial forms induces maps 
\[
  \wedge :  T_n^{p_1,q_1}(\Z)  \otimes T_n^{p_2,q_2}(\Z)  \to T_{n}^{p_1+p_2,q_1+q_2}(\Z)  \, . 
\]
We define a coboundary map  
\[
   d : T_{n}^{p,q}(\Z) \to T_{n}^{p+1,q}(\Z) 
\]
by $x_i \mapsto dx_i$, $dx_i \mapsto 0$ and extend it  to the whole of $T_n^{p,q}(\Z)$ by the Leibniz rule. 
Upon defining 
\[
   T_{n}^{p,q} := T_n^{p,q}(\Z)  \otimes \Q_q 
\]  
we obtain a filtered CGDA $T_n^{*}$ for each $n \geq 0$. Via pullback of forms, the maps $\delta_i$ and $\sigma_i$ 
induce maps of filtered CGDAs 
\begin{eqnarray*} 
   \partial_i : T_{n}^{*}  &  \to &  T_{n-1}^{*}  \\
      s_i : T_{n}^{*}         &  \to  &  T_{n+1}^{*}
\end{eqnarray*}
that satisfy  the simplicial identities. In other words, $(T_{n}^{*})_{n \in \N}$ is a simplicial 
filtered CGDA. 

Let $\mathcal{S}$ denote the category of simplicial sets. For any simplicial set $X \in \mathcal{S}$, 
we define a filtered CGDA 
\[
    T^{*,q}(X) := \Mor_{\mathcal S} (X, T^{*,q}) \, , 
\]
which is called the {\em Cenkl-Porter complex of polynomial forms on $X$}. Note that 
by definition $T^{*,0}(X)$ is concentrated in degree $0$ and can be identified with ${\rm Map}(\pi_0(X), \Z)$. 
Conversely, if 
$A \in \mathcal{A}$ is a filtered CGDA, its {\em simplicial realization} is the simplicial set
$\|A\|$ defined by
\[
    \| A \|_n  := \Mor_{\mathcal A} (A, T^{*}_n) 
\]
for $n \geq 0$. Note that corresponding constructions exist in rational homotopy theory, see \cite{BG}. 

For a proof of the following almost tautological statement, see \cite[8.1]{BG}. 

\begin{prop} \label{adjoint} Let $X$ be a simplicial set and $A$ be a filtered CGDA.  Then there is 
a canonical bijection 
\[ 
      \Mor_{\mathcal S} ( X, \|A\|) \approx \Mor_{\mathcal A}(A, T(X)) 
\]
which is natural in $A$ and $X$. In other words  the simplicial realization $\|-\|$ and 
the Cenkl-Porter functor $T(-)$ define a pair of right adjoint contravariant functors. 
\end{prop} 

Let us  fix the notation $\Psi_A : A \mapsto T(\|A\|)$ for the unit defined by this adjunction.  

We will now compare  $H^*(T^{*,q}(X))$ with the simplicial cohomology of $X$. This can be done most 
in the  axiomatic setting described in \cite{Ca}. We call a simplicial 
cochain complex  $(A^*_n)_{n \in \N}$ over a commutative ring $R$ with identity a {\em cohomology theory}, 
if the sequence of simplicial $R$-modules 
\[ 
    A^0 \stackrel{d}{\to}  A^1 \stackrel{d}{\to} A^2 \to \ldots 
\]
is exact, if the kernel $ZA^0$ of the differential $d : A^0 \to A^1$ is simplicially trivial, i.e.~all face and degeneracy maps are isomorphisms,  
and if for each $p \geq 0$, the simplicial  abelian group $A^p$ is contractible, cf.~\cite[Section 2]{Ca}. 
Concerning 
the latter property, recall \cite[Theorem 22.1]{May} that for a simplicial abelian group $(G_n)_{n \in \N}$, 
the $p$-th homotopy group $\pi_p(G)$ can 
be computed as the $p$-th homology of the chain complex 
\[
    \ldots \to G_n \to G_{n-1} \to G_{n-2} \to \ldots 
\]
with differentials equal to the alternating sum of the face operators of $G$. 
Given a cohomology theory $A^*$, we call the $R$-module $R(A) := (ZA^0)_0$ the {\em coefficients} of $A$. Furthermore 
we get an $R$-cochain complex 
\[ 
   A^*(X) := \Mor_{\mathcal S}(X, A^*)
\]
whose cohomology is naturally isomorphic to $H^*(X;R(A))$, the simplicial cohomology of $X$ with coefficients in $R(A)$, see \cite[Th\'eor\`eme 2.1]{Ca}.

\begin{prop} For each $q \geq 1$, the simplicial CGDA $T^{*,q}$ defines a cohomology theory with coefficients $\Q_q$. 
\end{prop} 

\begin{proof} The sequence 
\begin{equation*} 
     T^{0,q} \stackrel{d}{\to} T^{1,q} \stackrel{d}{\to} \ldots 
\end{equation*}
of simplicial $\Q_q$-modules is exact by the cubical 
Poincar\'e lemma \cite[Proposition 3.5]{CP}. The contractibility of  the simplicial sets $T^{p,q}$ 
for all $p \geq 0$ follows from the cubical extension lemma 
\cite[Lemma 3.1]{CP}.  That the kernel of the differential $T^{0,q} \to T^{1,q}$
is simplicially trivial is obvious (it consists of constant monomials). 
\end{proof} 
 
It is exactly at this point, where division  
by numbers smaller than or equal to $q$ is required: for $2 \leq k \leq q$ the closed form $t^{k-1} dt$ lives in filtration level $q$ and 
is the coboundary of $t^k/k$. Furthermore, notice that $T^{*,0}$ is not a cohomology theory. 
Indeed $T^{0,0} \to T^{1,0} \to T^{2,0} \to \ldots$ is exact and 
$T^{0,0}$ is simplicially trivial, but $T^{0,0}$ is not connected. 

The last  proposition implies that if $q \geq 1$ and $M$ is a $\Q_q$-module, 
then $M \otimes T^{*,q}$ is 
a cohomology theory with coefficients $M$. This follows from the universal coefficient theorem and 
the fact that $T^{*,q}$ consists of torsion free and hence flat $\Q_q$-modules, cf. \cite[Theorem 5.3.14]{Sp}.

\begin{cor} \label{kokett} Let $q \geq 1$ and let $M$ be a $\Q_q$-module. Then, for any $p \geq 0$, the simplicial $\Q_q$-module  $M \otimes T^{p,q}$ is 
contractible and for any $p \geq 1$, the simplicial $\Q_q$-module  $M \otimes ZT^{p,q}$   is an Eilenberg-MacLane 
complex of type $(M,p)$. 
\end{cor} 

\begin{proof} Both statements are implied by the fact that $M \otimes T^{*,q}$ is a cohomology theory. 
The first statement is then immediate. The second statement is proved by an inductive argument as in 
\cite[Th\'eor\`eme 1]{Ca} combined 
with the equation $M \otimes ZT^{p,q} =  \ker (\id \otimes d) : M \otimes T^{p,q} \to M \otimes T^{p+1,q}$ which 
follows from the universal coefficient theorem 
applied to the cochain complex of flat $\Q_q$-modules $0 \to  T^{p,q} \to T^{p+1,q} \to T^{p+2,q} \to \ldots$ 
\end{proof}

In order to study the relation of $H^*(T^{*,q}(X))$ and $H^*(X ; \Q_q)$, we first introduce for each $q \geq 1$ the cohomology 
theory $(C_n^{*,q})_{n \in \N}$, cf.~\cite[Section 3]{Ca}, with coefficients $\Q_q$, 
where $C^{*,q}_n$ is the simplicial cochain complex $C^*(\Delta[n] ; \Q_q)$ of the simplicial $n$-simplex. 
For each simplicial set $X$, we will identify $C^{*,q}(X)$ with $C^*(X;\Q_q)$, the simplicial cochain complex of $X$ with 
coefficients in $\Q_q$. 

For $q \geq 1$  we introduce a third simplicial  
cochain complex  $( (T \botimes C)^{*,q}_n)_{n \in \N}$ over $\Q_q$ with 
$(T \botimes C)^{*,q}_n := T^{*,q}_n \otimes  C^{*,q}_n$ and with face and degeneracy maps which are the tensor products 
of the corresponding maps on $T^{*,q}_n$ and $C^{*,q}_n$.  
The K\"unneth formula  \cite[Corollary 5.3.4]{Sp} 
and the Eilenberg-Zilber theorem \cite[Theorem 8.1]{MacLane} imply that 
$(T \botimes C)^{*,q}$ is again a cohomology theory with coefficients $\Q_q$. 
We have canonical maps $T^{*,q} \to  (T \botimes C)^{*,q}$ 
and $C^{*,q} \to (T \botimes C)^{*,q}$ of simplicial cochain complexes given by 
\begin{eqnarray*} 
    T_n^{*,q} &   = & T_n^{*,q} \otimes ( \Q_q \cdot 1) \hookrightarrow T_n^{*,q} \otimes C_n^{*,q}  \\
    C_n^{*,q} & = & (\Q_q \cdot 1) \otimes C_n^{*,q} \hookrightarrow T_n^{*,q} \otimes C_n^{*,q}\, .
\end{eqnarray*}
These induce isomorphisms of coefficients and hence \cite[Proposition 2]{Ca} we get induced isomorphisms
\begin{eqnarray*}
    H^*(T^{*,q}(X))  & \cong & H^*((T \botimes C)^{*,q}(X)) \\
    H^*(C^{*,q}(X)) & \cong & H^*((T \botimes C)^{*,q}(X)) \, , 
\end{eqnarray*}
which are natural in $X$, for each $q \geq 1$. Now we observe that exterior multiplication of cubical forms and 
the cup product of simplicial cochains define maps of simplicial cochain complexes
\begin{eqnarray*} 
      T_n^{*,q_1} \otimes T_n^{*,q_2} & \to &       T_n^{*, q_1 + q_2}  \\
      C_n^{*,q_1} \otimes C_n^{*,q_2} & \to &      C_n^{*,q_1 + q_2}  \\
      (T \botimes C)_n^{*,q_1} \otimes (T \botimes C)_n^{*,q_2} & \to & (T \botimes C)_n^{*,q_1+ q_2} 
 \end{eqnarray*}
which are compatible with the simplicial cochain maps 
$T_n^{*,q} \to (T \botimes  C)_n^{*,q}$ and $C_n^{*,q} \to (T \botimes C)_n^{*,q}$ considered before. We hence obtain 

\begin{thm}[Cenkl-Porter theorem \cite{CP}] \label{cenkl}  For $q \geq 1$ 
there are isomorphisms of $\Q_q$-modules  
\[
     H^*(T^{*,q}(X)) \cong H^*(X;\Q_q)  
\]
which are natural in $X$. These are compatible with the multiplicative structures in the sense that 
\[
 \xymatrix{
   H^{p_1}(T^{*, q_1}(X)) \otimes H^{p_2}(T^{*, q_2 }(X)) \ar[r]^-{\wedge}
   \ar[d]^{\cong} & 
   H^{p_1+p_2}(T^{*, q_1 + q_2 }(X)) \ar[d]^{\cong} \\
         H^{p_1}(X;\Q_{q-1}) \otimes H^{p_2}(X;\Q_{q_2}) \ar[r]^-{\cup}    &    H^{p_1 +p_2}(X;\Q_{q_1+q_2}) 
 }
\]
commutes for all $q_1, q_2 \geq 1$ and multiplicative units are mapped to each other. 
\end{thm} 

Recall that the isomorphism $H^*(T^{*,q}(X)) \cong H^*(X ; \Q_q)$ is induced by the zig-zag sequence of filtered graded differential 
algebras 
\begin{equation} \label{zigzag}
 T^{*}(X)  \rightarrow  (T \botimes C)^{*}(X)   \leftarrow    C^{*}(X)  
\end{equation}
where both arrows induce isomorphisms in cohomology.  In contrast, 
the proof of Theorem \ref{cenkl} in \cite{CP} is 
based on a single map  of filtered cochain complexes   
\[
 \int:   T^{*,q}(X) \to C^*(X, \Q_q)   
\]
induced by integration of forms. One can indeed show that $\int$ induces an isomorphism of 
cohomology groups which is compatible with multiplicative structures. However,  
$\int$ is not multiplicative on the cochain level.
Using the uniqueness of natural transformations of 
cohomology theories \cite[Proposition 2]{Ca} one can show that the isomorphisms induced by $\int$ and the one in our Theorem \ref{cenkl} coincide. 

Concerning filtration level zero, we note the existence of a canonical isomorphism 
\[
   H^*(T^{*,0}(X)) \cong H^0(X;\Z) \, . 
\]
Using this the multiplicativity property in Theorem \ref{cenkl} holds for $q_1, q_2  \geq 0$. 

As an addendum to the Cenkl-Porter theorem, we remark that for any $q \geq 1$ and 
for any prime $l > q$, we have isomorphisms 
\[
    H^*(T^{*,q}(X) \otimes \F_l) \cong H^*(X;\F_l)  \, , 
\]
which are multiplicative for filtration levels $q_1, q_2 \geq 1$ with $l > q_1 + q_2$.  This is shown by  
tensoring the cochain complexes appearing in (\ref{zigzag}) with $\F_l$ and applying the universal coefficient theorem. 

We finish this section with a discussion of higher Massey products. Let $l$ be an odd prime. It 
is well known \cite{Kraines} that the $l$-fold Massey product $H^1(B\Z/l; \F_l) \to H^2(B \Z/l; \F_l)$ 
is nonzero. Because squares of elements of degree one in commutative cochain algebras 
are equal to zero, one might wonder if this nonvanishing result does not contradict the Cenkl-Porter
theorem and in particular the addendum stated above. The following discussion will 
clarify this point.  

First we recall that Massey products 
are defined on singular cochain algebras $C^*(X;R)$ of topological spaces $X$, respectively 
simplicial cochain complexes of simplicial sets $X$, where 
$R$ is some commutative ring with unit. Because our proof of Theorem \ref{cenkl} is 
based on (\ref{zigzag}) in which all maps are multiplicative on the cochain level, we can use the filtered 
cochain algebra $T^{*,*}(X)$ or - in view of the addendum 
to the Cenkl-Porter theorem - the complexes $T^{*,q}(X) \otimes \F_l$ to construct higher Massey products 
in $H^*(X ; \F_l)$. But now the filtration structure of $T^{*,*}(X)$ has to be considered. 
More precisely, when constructing an $l$-fold Massey product of a cochain $c \in T^{1,1}(X) \otimes \F_l$, 
we have to perform $l-1$ multiplications (cf. the construction of Massey products in \cite{Kraines}) and 
hence arrive at elements in $T^{*,l}(X) \otimes \F_l$. But this group is zero, because $l$ is invertible in 
$T^{*,l}(X)$. 

This reasoning indicates that it is exactly the filtration structure together with its 
divisibility properties that allows us to generalize the construction of the Sullivan-de Rham algebra 
to an non-rational setting.

\section{Tame Hirsch Lemma} \label{Hirsch} 

In this section we shall explain a technique that will enable us to 
replace the Cenkl-Porter complex by a smaller filtered CGDA 
which reflects the topological structure of the 
underlying simplicial set in a similar manner as the minimal model does 
in rational homotopy theory.

Of fundamental importance is  the {\em filtration function}           
$\overline{\cdot} : \N_+ \to \N_+$ defined by 
\[
          \overline{1} = 1  \, , \hspace{0.5cm}  \overline{t} = 3 (t-1)  {\rm~for~}t\geq 2 \, . 
\]
One easily checks the following {\em admissibility} property: 
let $t_1, \ldots, t_n \geq 1$ and $n \geq 2$  be natural numbers. If $n = 2$, let  
$t_1, t_2 \geq 2$. Then the following implication holds for each $t \in \N_+$: 
\[
 \sum_{i=1}^n  t_i \leq  t + 1 \Longrightarrow \sum_{i=1}^n \overline{t_i} \leq \overline{t} \, . 
\]
It is easy to see that there is no function 
$\N_+ \to \N_+$ with  this property without this additional assumption if  $n = 2$. Futhermore, 
the above filtration function is optimal in the sense that any other admissible filtration 
function is pointwise greater than or equal to the one described above. 

In the literature on tame homotopy theory (see e.g.~\cite{Dw,Sch}), homotopy, homology and cohomology 
groups are indexed by degrees $r + k$, where $(r-1)$ denotes the connectivity of the spaces under consideration, and 
filtration functions depend on the parameter $t = k$. In our case we have $r=1$, and in order to simplify the notation,
we use the definition above which expresses the filtration function in terms of the actual degree $t = r+k$. 

\begin{defn} \label{techtech} Let $f : A^* \to B^*$ be a map of filtered CGDAs. 
The map $f$ is called a  {\em primary weak equivalence}, if for all primes $l>q$ the 
map $f \otimes \id : A^{*,q} \otimes \F_l \to B^{*,q} \otimes \F_l$
is a $t$-equivalence for all  $t \geq 1$ and $q \geq \overline{t}$, i.e. the induced map in cohomology is an 
isomorphism in degrees less than or equal to $t$ and injective in degree $t+1$.   . 

Given a natural number $k \geq 1$, we 
say that a primary weak equivalence $f : A^* \to B^*$ satisfies {\em condition $k^+$}, if for all primes $l >q$ 
the map $f  \otimes \id : A^{*,q} \otimes \F_l \to B^{*,q} \otimes \F_l $ 
is a $(t+1)$-equivalence for  all $t \geq k$ and $q \geq \max \{\overline{t}, \overline{k} +1 \}$. 
\end{defn} 

Note the close relationship between filtration level and cohomological degree. The ``extra degree'' in condition $k^+$ 
is essential for the inductive step in the construction of a small filtered cochain 
model along a Postnikov decomposition of the underlying simplicial set, see the proof of Lemma \ref{kplus}.

\begin{defn} Let $(B,d_B)$ be a filtered CGDA, $(V, d_V)$ a filtered 
cochain complex and $\tau : V^* \to B^{*+1}$ a filtration preserving cochain map of degree $1$, i.e.~$\tau$ satisfies 
the equation $d_B \circ \tau = - \tau \circ d_V$. 
The {\em free extension of $B$ by $V$ with twisting $\tau$}, denoted $B \otimes_{\tau} \Lambda V$, is 
the filtered graded algebra  $B \otimes \Lambda V$ equipped 
with the unique differential  $d$ which satisfies  $d(b \otimes 1) = d_B (b) \otimes 1$ for all $b \in B$ 
and $d(1 \otimes v) = \tau(v) \otimes 1 + 1 \otimes d_V(v)$ for all $v \in V$ and which acts as a graded derivation. 
\end{defn} 

Of particular importance are free extensions by so called elementary complexes. 

\begin{defn} Let $(p,q)$ be a pair of natural numbers with $p,q \geq 1$.  An {\em elementary complex} 
of {\em type} $(p,q)$ is a filtered cochain complex $(V^{*},\eta)$ such that 
\[
   V^{*,q'} = \begin{cases} 
                V^{*,q} \otimes \Q_{q'} {\rm~for~} q' \geq  q  \\
                0                       {\rm~for~} 0 \leq q' < q  \, ,    \end{cases}  
\]
and $V^{*,q}$ is of the form 
\[
    \ldots \to 0 \to V^{p,q} \stackrel{\eta}{\to} V^{p+1,q} \to 0 \to \ldots 
\] 
where $V^{p,q}$ and $V^{p+1,q}$ are finitely generated free $\Q_q$-modules and 
the cokernel of $\eta$ is a torsion module. 
\end{defn} 

If $V$ is an elementary complex of type $(p,q)$, then we let
\[
    \eta': (V^{p+1,q})' := \Hom(V^{p+1,q}, \Q_q) \to (V^{p,q})'
\]
denote the dual of $\eta$. By assumption the map $\eta'$ is injective. We 
will  show that simplicial realizations of elementary complexes are Eilenberg-MacLane
complexes for the group $\coker \eta'$. To this end let us first transfer some elementary constructions 
involving cochain complexes to the filtered setting.  

If $f : (V,d_V) \to (W,d_W)$ is a map of filtered cochain complexes, we denote by $(C_f, d_{C_f})$ 
the {\em mapping cone} of $f$, i.e. $C_f^{p,q} := V^{p,q} \oplus W^{p-1,q}$, where we 
set $W^{-1,q} = 0$,  and 
$d_{C_f}(v,w):= (d_V(v), f(v) - d_W(w))$. Special cases are 
the {\em suspension} $\Sigma W$ of $W$, defined  as the mapping cone of $0 \to W$, and 
the {\em cone}  over $V$, defined as the mapping cone of $\id : V \to V$ and 
denoted by $\cone V$. Note that we have a canonical short exact sequence 
\[
    0 \to \Sigma W    \to  C_f \to  V  \to 0  
\]
of filtered cochain complexes.

For the notions appearing in the following proposition and its proof see \cite{May}. 

\begin{prop} \label{eilen} Let $(V,\eta)$ be an elementary complex of type $(p,q)$. Then $\| \Lambda V \|$ 
is a (not necessarily minimal) Eilenberg-MacLane complex of type $(\coker \eta',p)$. In particular,
it satisfies the Kan condition. 
\end{prop}

\begin{proof}  By definition the set of $n$-simplices of $\| \Lambda V \|$ is given by 
\[
   \| \Lambda V \|_n = \Mor_{\mathcal A} (\Lambda V, T_n) = \Mor_{\mathcal C}(V, T_n) \, , 
\]
where $\mathcal{C}$ denotes the category of filtered cochain complexes. We conclude 
that $\| \Lambda V \|$ is a simplicial abelian group, thus fulfilling  the Kan condition, 
that fits into a pullback square of simplicial abelian groups 
\[
   \xymatrix{
        \| \Lambda V \| \ar[r] \ar[d] &   (V^{p,q})' \otimes T^{p,q} \ar[d]^{\id \otimes d} \\         
         (V^{p+1,q})' \otimes ZT^{p+1,q} \ar[r]^{\eta' \otimes \id} & (V^{p,q})' \otimes ZT^{p+1,q} 
   }   
\]
By Corollary \ref{kokett}, 
$(V^{p,q})' \otimes ZT^{p+1,q}$ is an Eilenberg-MacLane complex of type $((V^{p,q})',p+1)$ and 
$(V^{p,q})' \otimes T^{p,q}$ is contractible. Furthermore, the right 
hand vertical map is surjective with kernel $(V^{p,q})' \otimes ZT^{p,q}$. It is hence 
a principal Kan fibration with contractible total space and fibre an Eilenberg-MacLane complex 
of type $((V^{p,q})', p)$, cf.~\cite[Lemma 18.2]{May}. 
The long exact homotopy sequence for the induced fibration $\| \Lambda V \| \to (V^{p+1,q})' \otimes ZT^{p+1,q}$ 
shows that $\| \Lambda V \|$ is an Eilenberg-MacLane complex of type $(\coker \eta', p)$. 
\end{proof}

Given an elementary complex $(V, \eta)$ of type $(p,q)$ we will show that the simplicial realization functor transforms the sequence  
\[
       \Sigma V \to \cone V  \to V 
\]
into a fibre sequence 
\[
    \| \Lambda V \| \hookrightarrow \| \Lambda (\cone V) \| \to \| \Lambda (\Sigma V) \|  \, . 
\]
For this we note  that the map of simplicial abelian groups $\| \Lambda(\cone V) \| \to \| \Lambda(\Sigma V) \|$ 
is surjective and  hence a principal Kan fibration, whose kernel can be identified with the simplicial 
abelian group $\| \Lambda V \|$. The last point can be checked using  an explicit description of $\| \Lambda(\cone V) \|$ and 
$\| \Lambda(\Sigma V) \|$ similar to the one of $\| \Lambda V \|$ in the proof of Proposition \ref{eilen}. 
This explicit description also shows that the homotopy groups of $\| \Lambda(\cone V)\|$ are zero in all degrees
and hence that $\| \Lambda(\cone V) \|$ is contractible.  By Proposition \ref{eilen} the simplicial set $\| \Lambda(\Sigma V) \|$ is 
an Eilenberg-MacLane complex of type $(\coker \eta' , p+1)$. We conclude that 
the fibration 
\[
    \| \Lambda V \| \hookrightarrow \| \Lambda (\cone V) \| \to \| \Lambda (\Sigma V) \|
\] 
is fibre homotopy equivalent to the path-loop
fibration over an Eilenberg-MacLane complex of type $(\coker \eta',p+1)$. Fibrations of this kind can 
thus be used as  building block in the Postnikov decompositions of simplicial sets. 

\begin{prop} \label{adjointbasic} Let $V$ be an elementary complex of type $(p,q)$ with $V^{p,q} = \Q_q$
and $V^{p+1,q} = 0$. Then 
the unit $\Psi_{\Lambda V} : \Lambda V \to T(\| \Lambda V \|)$ induces an isomorphism 
$H^p((\Lambda V)^{*,q}) \to H^p(T^{*,q}( \| \Lambda V \| ) )$ of groups that can be 
canonically identified with $\Q_q$. 
\end{prop} 

\begin{proof} The generator of $H^{p}(T^{*,q}( ZT^{p,q})) = \Q_q$ is represented by 
the inclusion map $ZT^{p,q} \hookrightarrow T^{p,q}$ in   $\Mor_{\mathcal S}(ZT^{p,q}, T^{p,q})$. 
Now recall from the proof of Proposition \ref{eilen} that $\|\Lambda V\| = (V^{p,q})' \otimes ZT^{p,q}$.
We thus  obtain the unit  
\[
   \Psi_{\Lambda V} :   (\Lambda V)^{p,q} \to  \Mor_{\mathcal S}\big( (V^{p,q})' \otimes ZT^{p,q}, T^{p,q} \big) 
\]
and this sends $c  \in (\Lambda V)^{p,q} = \Q_q$ to the morphism of simplicial sets 
$(V^{p,q})' \otimes ZT^{p,q} \to T^{p,q}$, $\phi \otimes  z_n \mapsto \phi(c) \cdot z_n$. Together 
with the previous remark this implies the assertion. 
\end{proof} 

Now let $X$ be a simplicial set, let $V$ be an elementary complex and let 
\[
   f : X \to \| \Lambda(\Sigma V) \|
\]
be a simplicial map. We obtain a pull back square
\[
\xymatrix{
     E_f \ar[r] \ar[d]^{p_f} & \|\Lambda(\cone V) \| \ar[d] \\ 
      X  \ar[r]^-{f}  &  \|\Lambda(\Sigma V)  \| 
}  
\]
The adjoint of $f$ is a map of filtered CGDAs $\Lambda(\Sigma V) \to T(X)$   
and this is induced by a filtered cochain map $f^{\sharp} : V \to T(X)$ of degree $1$. We hence get a 
free extension $T(X) \otimes_{f^{\sharp}}  \Lambda V$ and the maps $T(p_f) : T(X) \to T(E_f)$ and 
the composition $V \hookrightarrow \cone V  \to T(E_f)$ induced by the adjoint of $E_f \to \| \Lambda (\cone V) \|$  
induce a filtered CGDA map 
\[
    \Gamma_f : T(X) \otimes_{f^{\sharp}} \Lambda V \to T(E_f) \, . 
\] 

The following result will be important in the next section. It is a combination of \cite[Theorem II.6.4.(i)]{So}  and 
\cite[Lemma II.8.0]{So}. A detailed proof will be given in Section \ref{append}.  

\begin{thm}[Tame Hirsch lemma] \label{tamehirsch}
Let $V$ be an elementary complex of 
type $(d, \overline{k})$, where $d \geq k \geq 1$,  
let $X$ be a simplicial set and let $f : X \to \|\Lambda(\Sigma V)\|$ be a simplicial map. Then the induced map 
\[
   \Gamma_f : T(X) \otimes_{f^{\sharp}} \Lambda V \to T(E_f) 
\]
is a primary weak equivalence satisfying condition $k^+$. 
\end{thm}

We conclude this section with an example which illustrates a fundamental  problem 
with elements of degree $1$ and justifies the specific form of Definition \ref{techtech}. 
Let $V$ be an elementary 
complex of type $(1,1)$, concentrated in degree $1$ with $V^{1,1} = \Z$. Then $\|\Lambda V\|$ is an 
Eilenberg-Mac Lane complex 
of type $(\Z,1)$ and for each $q \geq 1$, the unit $(\Lambda V)^{*,q} 
\to T^{*,q}(\| \Lambda V\|)$ induces isomorphisms 
in cohomology in all degrees. 
For a constant map  $f : \| \Lambda V \| \to \|\Lambda ( \Sigma V) \|$, we study the map   
\[
    \Gamma_f :  T(\|\Lambda V \|) \otimes \Lambda V \to T(\| \Lambda V \| \times \| \Lambda V \|) 
\]
that appears in the tame Hirsch lemma.  This map is not a primary $2$-equivalence in filtration level $1$,
because $H^2((  T ( \| \Lambda V\|) \otimes \Lambda V)^{*,1})=0$, but 
$H^2(T^{*,1}(\| \Lambda V \| \times \| \Lambda V \|)) = \Z$
 by the Cenkl-Porter theorem. This example shows that the 
extra degree in property $k^+$ does not occur below filtration level $\overline k+1$ in general.

\section{Small cochain models for $l$-tame $(\Z/l)^r$-spaces}  \label{small} 

Let $l$ be a prime, $G := (\Z/l)^r$ and $X$  a $G$-space. 
Under certain assumptions on $X$ we shall construct a small commutative graded 
differential algebra over $\F_l$ which is
free as a graded algebra and whose 
cohomology is  equal to $H^*(X_G; \F_l)$ up to a certain degree.   
In particular we will prove Proposition \ref{later} from Section \ref{cool}. 

From now on it will be  convenient to work exclusively in an $l$-local setting.  This means that all 
filtered cochain complexes and filtered cochain algebras appearing in Sections \ref{CenklPorter} and 
\ref{Hirsch} will be  tensored with $\Z_{(l)}$ in each filtration level. Since  tensoring with 
$\Z_{(l)}$ is  exact, we still obtain filtered objects. For example, 
if $V$ is a filtered cochain complex, then it is understood from now on that 
$V^{*,q}$ is a complex of $\Z_{(l)}$-modules for $0 \leq q < l$ and 
a complex of $\Q$-modules for $q \geq l$.  The Cenkl-Porter theorem now states that 
we get natural multiplicative isomorphisms 
\[
    H^*(T^{*,q}(X)) \cong \begin{cases} H^*(X;\Z_{(l)})   {\rm~for~}1 \leq q < l   \\
                                        H^*(X;\Q)         {\rm~for~}q \geq l \, ,  
                          \end{cases} 
\]
where filtration levels are additive under multiplication of elements.  
Of course, in the definition and discussion of primary weak equivalences satisfying condition $k^+$, one now restricts to the single 
prime $l$ and to filtration levels $q < l$. This is  underlined by the notion  {\em $l$-primary weak equivalence satisfying condition $k^+$}, 
which will be used in the sequel.  

\begin{defn} \label{tamedef} We call $X$ an {\em $l$-tame} $G$-space, if $X$ is path connected, 
$\pi_1(X)$ is abelian and for all $t \geq 1$ with $\overline{t} < l$ 
the following holds: 
\begin{enumerate}[1)]
  \item $\pi_t(X)$ is a finitely generated free $\Z_{(l)}$-module, 
  \item $\pi_1(X)$ acts trivially on $\pi_t(X)$ and
  \item the induced $G$-action on $\pi_t(X)$ is trivial. 
\end{enumerate} 
\end{defn}

\begin{lem} \label{lsimple} Let $X$ be an $l$-tame $G$-space. Then $\pi_1(X_G)$ is abelian and for each $t \geq 2$ with $\overline t < l$, 
$\pi_t(X_G)$ is a finitely 
generated free $\Z_{(l)}$-module on which $\pi_1(X_G)$ acts trivially. 
\end{lem} 

\begin{proof} We only show that $\pi_1(X_G)$ is abelian. 
The proofs of the other assertions are elementary and left to the reader. 

The short exact sequence $1 \to \pi_1(X) \to \pi_1(X_G) \to G  \to 1$ 
can be used to construct a map 
\[
   \phi :  G \times G \to \pi_1(X) \, , ~~ (g_1, g_2) \mapsto \overline g_1 \cdot \overline g_2 \cdot (\overline g_1)^{-1}
                                                    \cdot (\overline g_2)^{-1} \, . 
\]
Here $\overline g_1$ and $\overline g_2$ are lifts of $g_1, g_2 \in G = (\Z/l)^r$ to $\pi_1(X_G)$.
The map $\phi$ is well-defined because the conjugation action of  $G$ on $\pi_1(X)$ is trivial. Futhermore,  $\phi$ 
is $\Z/l$-bilinear. Because $G$ is a torsion group and $\pi_1(X)$ is a free 
$\Z_{(l)}$-module, $\phi$ must be the zero map and hence $\pi_1(X_G)$ is abelian. 
\end{proof}

Later  we will need the following result, which is important to study maps between 
cochain complexes over  $\Z_{(l)}$. 

\begin{lem} \label{univkoeff} Let $l$ be a prime and let $V^*$ and $W^*$ be (unfiltered) 
cochain complexes over $\Z_{(l)}$ that are torsion free $\Z_{(l)}$-modules in each degree. 
Furthermore assume that $H^*(V^*)$ and $H^*(W^*)$ are 
finitely generated $\Z_{(l)}$-modules in each degree. 

Now let $f : V^* \to W^*$ be a map of cochain complexes so that
$f \otimes \id: V^* \otimes \F_l \to W^* \otimes \F_l$ is a $t$-equivalence for some $t \geq 1$.  
Then $f^* : V^* \to W^*$ is itself a $t$-equivalence. 
\end{lem} 

\begin{proof} The short exact sequence 
\[
   0 \to \Sigma W^* \to C_f \to V^* \to 0 
\]
induces long exact sequences in cohomology 
\[
   \ldots \to H^i(C_f \otimes M) \to H^i(V \otimes M) \stackrel{f^*}{\to} H^i(W \otimes M) \to \ldots 
\]
for any $\Z_{(l)}$-module $M$. Our assumptions imply that 
$H^i(C_f \otimes \F_l) = 0$ for all $0 \leq i \leq t+1$. Also, $C_f$ consists of flat $\Z_{(l)}$-modules in each degree. 
Hence we can apply the universal coefficient theorem \cite[Theorem 5.3.14]{Sp}  to conclude that $H^i(C_f)\otimes \F_l = 0$  and 
hence that $H^i(C_f) = 0$ for 
$0 \leq i \leq t+1$. Here we use that  $H^*(C_f)$ is finitely generated over $\Z_{(l)}$ in 
each degree. From this the assertion follows.  
\end{proof}

From now on we work in a simplicial setting and assume that we are 
given a Kan complex $X$ equipped with a $G$-action. This can be achieved by passing to the 
simplicial set of singular simplices of a given $G$-space. 
As our model for the universal fibration $EG \to BG$ we take the universal covering 
$EG \to \| \Lambda V_0 \|$, where  $V_0$ is the elementary complex of type $(1,1)$ 
given on filtration level $1$ by the cochain complex 
\begin{equation} \label{vaunull}
   0 \to \Z_{(l)}^r \stackrel{c \mapsto l \cdot c}{\to} \Z_{(l)}^r \to 0 \to 0 \to \ldots \, , 
\end{equation} 
compare Proposition \ref{eilen}. In this way we realize the Borel fibration 
$X \hookrightarrow X_G \to BG$ as a Kan fibration.  

Now let $X$ be an $l$-tame $G$-complex. We will study the Postnikov decomposition of $X_G$ relative to $BG$. 
This leads to a commutative diagram 
\[
 \xymatrix{ 
     \ar@{.>}[d]            &       \ar@{.>}[d]                  &  \ar@{.>}[d]  \\
  X_2 \ar[r] \ar[d]^{p_2}  &      (X_G)_2 \ar[r] \ar[d]^{P_2}  &   BG \ar@{=}[d] \\
  X_1 \ar[r] \ar[d]^{p_1}  &      (X_G)_1 \ar[r]^{P_1} \ar[d]^{P_1}  &   BG \ar@{=}[d] \\
X_0 = \star \ar[r]     &    (X_G)_0 = BG     \ar@{=}[r] &   BG 
  }
\]
where for all $k \geq 1$ the complexes $X_k$ and $(X_G)_k$ 
are $k$-th stages of the Postnikov decompositions of $X$ and $X_G$, 
each row is a fibre sequence with fibre $X_k$ and the vertical maps $p_k$ and $P_k$ are Kan fibrations whose fibres 
are Eilenberg-MacLane  complexes of type $(\pi_k(X),k)$.  
By the assumption that $X$ is $l$-tame together with Lemma \ref{lsimple},
the maps $X_k \to X_{k-1}$ and $(X_G)_k \to (X_G)_{k-1}$ 
are simple fibrations, if $k\geq 1$ and $\overline{k} < l$, i.e the fundamental group of the base acts trivially on the $k$-th homotopy      
group of the fibre. Therefore, by Proposition \ref{eilen} and the subsequent remarks,  
we can assume that these maps fit into diagrams 
\[
  \xymatrix@=12pt{
       &      \| \Lambda(\cone V_k) \|  \ar@{=}[rr] \ar'[d][dd] &     &    \| \Lambda(\cone V_k) \| \ar[dd]    \\
   X_{k} \ar[ur] \ar[rr] \ar[dd]^<<<<{p_k}  &       &          (X_G)_{k} \ar[ur] \ar[dd]^<<<<{P_k}  \\
      &      \| \Lambda(\Sigma V_k)  \|  \ar@{=}'[r][rr] &  & \| \Lambda(\Sigma V_k) \|  \\
   X_{k-1}     \ar[rr] \ar[ur]_<<<<<<<<{f_{k-1}}    &       &          (X_G)_{k-1} \ar[ur]_<<<<<<<<{F_{k-1}}  
  }
\]
where $(V_k,\eta)$ is an elementary complex of type $(k,\overline{k})$ 
which is concentrated in degree $k$ with $V^{k, \overline k} = \pi_k(X)$.

For $k = 0$ and for $k \geq 1$, $\overline{k} < l$, we will define filtered CGDAs 
$M_k$ and $E_k$ together with filtered CGDA maps 
\begin{eqnarray*}
    \phi_k : E_{k} & \to & M_{k}  \\
    \psi_k : M_k   & \to & T(X_k) \\
    \Psi_k : E_k   & \to & T((X_G)_k)   
\end{eqnarray*} 
fitting into commutative diagrams   
\[
  \xymatrix@=25pt{
       &      T(X_k)                  &     &    T((X_G)_k)   \ar[ll]         \\
   M_{k} \ar[ur]^-{\psi_k}                    &       &             E_k  \ar[ur]^-{\Psi_k} \ar[ll]_<<<<<<<<<<{\phi_k}    \\
      &       T(X_{k-1})   \ar'[u]^<<<<<<{T(p_k)}[uu] &     &    T((X_G)_{k-1}) \ar'[l][ll] \ar[uu]^<<<<<<{T(P_k)}  \\
   M_{k-1}       \ar[ur]^-{\psi_{k-1}} \ar[uu]   &       &        E_{k-1} \ar[ll]_<<<<<<<<{\phi_{k-1}} \ar[ur]^-{\Psi_{k-1}} \ar[uu]    
  } 
\]
Furthermore, for $k \geq 1$ the maps  $\psi_k : M_k \to T(X_k)$ and $\Psi_k : E_k \to T((X_G)_k)$ 
will be $l$-primary weak equivalences satisfying condition $k^+$, the maps $\psi_0$ and $\Psi_0$ will be $l$-primary weak equivalences 
satisfying condition $1^+$ and  the vertical maps in the front square will be inclusions which can be used to express $M_k$ and $E_k$ 
as free extensions of $M_{k-1}$ and $E_{k-1}$ making the following diagram commutative:
\[
  \xymatrix@C=50pt@R=30pt{
              M_k         &       E_k \ar[l]_{\phi_k} \\
    M_{k-1} \otimes_{\phi_{k-1} \circ \tau_k} \Lambda V_{k} \ar@{=}[u]  &   E_{k-1} \otimes_{\tau_k} \Lambda V_k 
    \ar@{=}[u] \ar[l]_-{\phi_{k-1} \otimes \id}
  } 
\]
We define these objects by induction on $k$. Let 
$M_0^{*,q} = \Z_{(l)}$, concentrated in degree $0$,   for $0 \leq q < l$, and let $M_0^{*,q} = \Q$,
concentrated in degree $0$,  
for $q \geq l$. Furthermore, define  $E_0 := \Lambda V_0$, where the cochain complex $V_0$ was 
defined in (\ref{vaunull}). The maps  
\begin{eqnarray*}
      \phi_0 : E_0 & \to & M_0  \\
      \psi_0 : M_0 & \to & T(X_0) = T(\star)  \\
      \Psi_0 : E_0 & \to & T((X_G)_0) = T(\|V_0\|) 
\end{eqnarray*}
are the obvious ones: $\phi_0$ is an isomorphism in degree $0$ and $\Psi_0$ is the unit. 
The map $\Psi_0$ is an $l$-primary weak equivalence satisfying condition $1^+$ by 
Proposition \ref{hui}. Now let $k \geq 1$ and $\overline{k} < l$ and assume 
that the above objects have been constructed for $k-1$. 

\begin{lem} \label{kplus} 
For $k \geq 1$ and $\overline k < l$  the map $\Psi_{k-1}$ induces an  isomorphism 
\[
    \Psi_{k-1}^* : H^{k+1}(E_{k-1}^{*,\overline k})  \cong   H^{k+1}\big( T^{*, \overline k }((X_G)_{k-1}) \big) \, . 
\]
\end{lem} 

\begin{proof} At first we notice that the cochain complexes $E_{k-1}^{*,\overline k}$ and 
$T^{*, \overline k} ((X_G)_{k-1})$ consist of torsion  free $\Z_{(l)}$-modules in each degree.  
Moreover,  $H^*(E_{k-1}^{*,\overline k})$  and $H^*(T^{*,\overline k} ((X_G)_{k-1}))$ 
are finitely generated over $\Z_{(l)}$ in each degree. This holds by the Cenkl-Porter theorem 
and because $\pi_*((X_G)_{k-1})$ is finitely generated over $\Z_{(l)}$  in each degree. Remember that $X$ is $l$-tame 
and that $\overline k < l$.  

Now, for $k \geq 2$ the assertion follows by Lemma \ref{univkoeff}, 
because $\Psi_{k-1}$ is an $l$-primary weak equivalence satisfying condition $(k-1)^+$ and hence an $l$-primary $(k+1)$-equivalence 
in filtration level $\overline k$. 

For $k=1$ we get $E_{k-1} = \Lambda V_0$ and  the composition 
\[
   (\Lambda V_0)^{*,\overline 1} \otimes \F_l \to T^{*,\overline 1}(\| \Lambda V_0 \|) \otimes \F_l 
\]
is a $1$-equivalence by Proposition \ref{hui}.  
An application of Lemma \ref{univkoeff} shows that the induced map 
\[
   H^{2}( (\Lambda V_0)^{*,\overline 1})    \to  H^2 ( T^{*,\overline 1}(\| \Lambda V_0 \|))   
\]
is injective. It is hence an isomorphism as both  source and target 
are isomorphic to $(\Z/l)^r$.  
\end{proof} 

Notice that in view of the remarks at the end of Section \ref{Hirsch} 
it is somewhat unexpected that we still get an isomorphism for $k=1$.

Because $H^{*}(\Sigma V_k)$ is a free $\Z_{(l)}$-module concentrated in degree $k+1$, Lemma \ref{kplus} shows that the map 
\[
          H^*((\Sigma V_k)^{*,\overline k}) \to H^*(\Lambda (\Sigma V_k)^{*,\overline k}) \to H^*(T^{*,\overline k}((X_G)_{k-1})) 
\]
induced by the adjoint of  $F_{k-1}$ can be written as the composition of $\Psi_{k-1}^*$ and a map 
\[
  H^*((\Sigma V_k)^{*,\overline k}) \to H^*(E_{k-1}^{*, \overline k}) \, . 
\]
This map is 
induced by a map of filtered cochain complexes
\[
   \tau_k :  \Sigma V_k^{*} \to E_{k-1}^{*} 
\]
and defines the free extensions $E_k = E_{k-1} \otimes_{\tau_k} \Lambda V_k$ and 
$M_k = M_{k-1} \otimes_{\phi_{k-1} \circ \tau_k} \Lambda V_k$ as well as the map $\phi_k := (\phi_{k-1} \otimes \id) :  E_k \to M_k$. 

It remains to construct the $l$-primary weak equivalences satisfying condition $k^+$ 
\begin{eqnarray*}
   \psi_k : M_k  & \to &  T(X_k) \\
   \Psi_k : E_k  & \to  &  T((X_G)_k) \, . 
\end{eqnarray*}
Let 
\[  
   F_{k-1}^{\sharp}:     \Sigma V_k  \to T((X_G)_{k-1}) 
\]
be the restriction of the adjoint of $F_{k-1}$.  

The composition
\[
           \Psi_{k-1} \circ  \tau_k   : \Sigma V_k^{*,\overline k} \to T^{*, \overline k} ((X_G)_{k-1} )
\]
induces the same map in cohomology as the map $F_{k-1}^{\sharp}$. We denote by $\overline{F_{k-1}}$ the 
adjoint of $\Psi_{k-1} \circ \Lambda(\tau_k) : \Lambda (\Sigma V_k) \to E_{k-1} \to T((X_G)_{k-1})$
and conclude that the maps $F_{k-1}$ and $\overline{F_{k-1}}$ are homotopic as maps 
of simplicial sets $(X_G)_{k-1} \to \| \Lambda ( \Sigma V_k )\|$. This follows, because
 $\| \Lambda ( \Sigma V_k )\|$ is an Eilenberg-MacLane complex of type $(\pi_k(X), k+1)$ and 
because after applying the functor $T^{*, \overline k}$ 
and precomposing with the unit 
\[
   \Psi_{\Lambda(\Sigma V_k)}|_{\Sigma V_k^{*,\overline k}} : \Sigma V_k^{*,\overline k} \to T^{*,\overline k}(\| \Lambda (\Sigma V_k) \|) \, , 
\]
we precisely  obtain the maps $F_{k-1}^{\sharp}$ and 
$\Psi_{k-1} \circ \tau_k$ considered at the beginning of this paragraph, see also Proposition \ref{adjointbasic}.  

By precomposing with the inclusion $X_{k-1} \to (X_G)_{k-1}$ we also 
obtain a homotopy between $f_{k-1}$ and $\overline{f_{k-1}}$ 
where the last map is adjoint to $\psi_{k-1} \circ \phi_{k-1} \circ \Lambda(\tau_k)$.  
It follows that there are homotopy equivalences 
\begin{eqnarray*} 
   \xi :   X_k  & \simeq   & E_{\overline{f_{k-1}}} \\
   \Xi :   (X_G)_k  & \simeq  & E_{\overline{F_{k-1}}}  
\end{eqnarray*}
making the diagram 
\[
   \xymatrix{
             X_k \ar[rrr] \ar[dd]_{\xi} \ar[dr]^{p_k}  & & & (X_G)_k \ar[dd]_{\Xi} \ar[dl]_{P_k} \\
                                   &  X_{k-1} \ar[r] & (X_G)_{k-1} \\
             E_{\overline{f_{k-1}}} \ar[rrr] \ar[ur] & & & E_{\overline{F_{k-1}}} \ar[ul] 
   }
\]
commuative.   
We now define $\Psi_k$ and $\psi_k$ as the compositions
\begin{eqnarray*}
  &&  \Psi_k :  E_{k-1} \otimes_{\tau_k} \Lambda V_k \stackrel{\Psi_{k-1} \otimes \id}{\longrightarrow} 
   T((X_G)_{k-1}) \otimes_{\Psi_{k-1} \circ \tau_k} \Lambda V_k 
          \stackrel{\Gamma_{\overline{F_{k-1}}}}{\to} T(E_{\overline F_{k-1}}) \stackrel{T(\Xi)}{\longrightarrow} T((X_G)_k)   \\
  &&  \psi_k : M_{k-1} \otimes_{\phi_{k-1} \circ \tau_k} \Lambda V_k \stackrel{\psi_{k-1} \otimes \id}{\longrightarrow} \, . 
         T(X_{k-1}) \otimes_{\psi_{k-1} \circ \phi_{k-1} 
\circ \tau_k} \Lambda V_k \stackrel{\Gamma_{\overline{f_{k-1}}}}{\longrightarrow} T(E_{\overline f_{k-1}}) 
         \stackrel{T(\xi)}{\to} T(X_k)    
\end{eqnarray*}
By the tame Hirsch lemma and Lemma \ref{aufwand}, these maps are $l$-primary weak equivalences satisfying condition $k^+$, since 
a primary weak equivalence satisfying condition $(k-1)^+$ is also a primary weak equivalence satisfying condition $k^+$ and 
$(X_G)_{k-1}$ is connected by assumption. 
This finishes the construction of $\Psi_{k}$ and $\psi_k$.

We now formulate the main result of this section.

\begin{thm} \label{smallcochain} Let $X$ be an $l$-tame $G$-space  and let $k \geq 1$ be a number 
with $\overline{k} < l$. Then there are (unfiltered) commutative graded differential algebras $(E^*,d_E)$ 
and $(M^*,d_M)$ over $\F_l$ with the following properties: 
\begin{enumerate}[1)]
   \item as graded algebras $E^* = \F_l[t_1, \ldots, t_r] \otimes \Lambda(s_1, \ldots, s_r) \otimes M^*$
         where each $t_{i}$ is of degree $2$ and each $s_{i}$ is of degree $1$.  
   \item the differential $d_E$ is zero on $\F_l[t_1, \ldots, t_r] \otimes \Lambda(s_1, \ldots, s_r) \otimes 1$ 
         and the map $E^* \to M^*$, $t_i , s_i \mapsto 0$, is a cochain map.  
   \item $M^*$ is free as a graded algebra with generators in degrees $1, \ldots , k$. As generators in degree 
                  $i \in \{1, \ldots, k\}$ we can take the elements of a basis of the free  $\Z_{(l)}$-module $\pi_i(X)$. 
   \item the cohomology algebra of $E^*$ is multiplicatively isomorphic to $H^*(X_G; \F_l)$ up to degree $k$. 
   \item assume that $H^{*}(X_G ; \F_l)$ vanishes in degrees $k$ and $k+1$. Then each monomial in 
             $t_1, \ldots, t_r$ 
             of cohomological degree at least $k$ represents the zero class in $H^*(E)$.  
\end{enumerate} 
\end{thm} 

\begin{proof} After $k$ steps the above process leads to  a commutative diagram 
\[
    \xymatrix{
                 E_k \ar[r]^{\phi_k}  \ar[d]^{\Psi_k}  &  M_k \ar[d]^{\psi_k} \\
                     T((X_G)_k) \ar[r]  &  T(X_k) 
    }
\]
where the vertical arrows are $l$-primary weak equivalences satisfying condition $k^+$. We can restrict this diagram 
to filtration level $\overline{k}$ to obtain a commutative diagram of $\F_l$-cochain complexes 
\[  
\xymatrix{
        E_k^{*,\overline k} \otimes \F_l \ar[r] \ar[d] & M_k^{*,\overline k} \otimes \F_l 
 \ar[d] \\
     T^{*, \overline k}((X_G)_k)  \otimes \F_l \ar[r] & T^{*, \overline k}(X_k) \otimes \F_l \, .  
}
\]
The vertical arrows in this diagram are $k$-equivalences. By its inductive construction 
the cochain complex $E_k^{*, \overline k}$ can be written as
\[
    E_k^{*, \overline k} = \big(  ((\Lambda V_0 \otimes_{\tau_1} \Lambda V_1) \otimes \ldots ) \otimes_{\tau_k} \Lambda V_k \big)^{*, \overline k} 
\]
where $V_0$ is an elementary complex of type $(1,1)$ and $V_t$ is an elementary complex of type $(t, \overline t)$ 
which is concentrated in degree $t$ for $t \in \{1, \ldots, k\}$.

We set 
\[
 W_0 := \Lambda ( V_0^{*,1} \otimes \F_l) \, , \hspace{.5cm}  W_t  :=   \Lambda (V_t^{*,\overline t} \otimes \F_l) 
{\rm~~for~}1 \leq t \leq k \, ,  
\]
regarded as unfiltered $\F_l$-cochain algebras. We hence get an unfiltered $\F_l$-cochain algebra 
\[
    E^* = ( ( W_0 \otimes_{\tau_1} W_1 ) \otimes \ldots  ) \otimes_{\tau_k} W_k \, . 
\]
Note that $W_0 = \F_l[t_1, \ldots, t_r] \otimes \Lambda(s_1, \ldots, s_r)$ with trivial 
differential. We denote the evaluation  $t_i, s_i \mapsto 0$ by $\phi$ and obtain an  $\F_l$-cochain 
algebra 
\[
   M^* = ( W_1 \otimes_{\phi \circ \tau_2} W_2 \otimes \ldots ) \otimes_{\phi \circ \tau_k} W_k \, . 
\]
It follows that $M^*$ is of the form described in Theorem \ref{smallcochain}, that 
\[
    E^* = \F_l[t_1, \ldots, t_r] \otimes \Lambda (s_1, \ldots, s_r) \otimes M^* 
\]
as graded algebras and that $E^* \to M^*$, $t_i ,s_i \mapsto 0$, is a cochain map. 

It remains to study the relation between $H^*(E)$ and $H^*(X_G;\F_l)$ as stated in parts 4)  and 5) 
of Theorem \ref{smallcochain}. Here we have to keep in 
mind that the product 
\[
   (E_k^{p_1,q_1} \otimes \F_l) \otimes (E_k^{p_2,q_2} \otimes \F_l) \to E_k^{p_1 + p_2, q_1+q_2} \otimes \F_l
\]
is the zero map, if $q_1 + q_2 \geq l$. Therefore we have to be somewhat careful when comparing the 
multiplications in the unfiltered algebra $E^*$ with those in the filtered algebra $E_k^{*}$. 

Let $B_0\subset V_0^{*,1}$ and $B_t \subset V_t^{t, \overline t}$  be  homogenous $\Z_{(l)}$-module
bases for $1 \leq t \leq k$. We consider the set $S$ of monomials of homological degree 
at most $k$ in elements of $B_0 \cup \ldots \cup B_k$. 
The admissibility of the filtration function implies that 
$S \subset E_k^{*,\overline k}$. Let $Z^* \subset E_k^{*,\overline k}$ 
be the $\Z_{(l)}$-subcomplex which is generated by the monomials in $S$ and their coboundaries. 

We obtain inclusions $Z^* \otimes \F_l \to E^*$ and $Z^* \otimes \F_l \to E_k^{*, \overline k} \otimes \F_l$ 
of cochain complexes both of which are $k$-equivalences. This 
leads to an additive isomorphism 
\[
    \alpha : H^{\leq k}(E^*) \cong H^{\leq k}(E_{k}^{*, \overline k} \otimes \F_l) \, . 
\]
Consider the composition
\[
 \beta : H^*(E_k^{*,\overline{k}} \otimes \F_l) \to H^*(T^{*,\overline{k}} ((X_G)_k) \otimes \F_l ) \stackrel{\cong}{\to}
   H^*((X_G)_k ; \F_l) \to H^*(X_G ; \F_l)  
\]
where the second map is the Cenkl-Porter map and the last map is induced 
by the canonical map  $X_G \to (X_G)_k$ to the $k$-th Postnikov piece. As a composition of 
$k$-equivalences the map $\beta$  is itself a $k$-equivalenc.

In summary, we obtain an additive isomorphism 
\[
 \beta \circ \alpha :   H^{\leq k}(E^*) \cong H^{\leq k}(E_{k}^{*, \overline k} \otimes \F_l) \cong  H^{\leq k}(X_G;\F_l) \, . 
\]
In order to study the multiplicative behaviour of $\beta \circ \alpha$, let 
$c_1 \in H^{p_1}(E)$ and $c_2 \in H^{p_2}(E)$, where $1 \leq p_1, p_2 \leq k$ and $p_1 + p_2 \leq k$.
Then, by the construction of $E_k^{*, \overline k}$, the elements $\alpha(c_1)$ and $\alpha(c_2)$ 
can be represented by cocycles 
$\gamma_1 \in E_{k}^{p_1, \overline{p_1}} \otimes \F_l$ 
and $\gamma_2 \in  E_{k}^{p_2, \overline{p_2}} \otimes \F_l$ 
and hence $\alpha(c_1 \cdot c_2)$ can be represented by the 
cocycle $\gamma_1 \cdot \gamma_2 \in E_{k}^{p_1+ p_2, \overline{p_1} + \overline{p_2}} \otimes \F_l 
\subset E_{k}^{p_1 + p_2, \overline{p_1 + p_2}} \otimes \F_l \subset E_{k}^{p_1+p_2, \overline{k}} \otimes \F_l$. Here 
we have used the assumption $p_1 + p_2 \leq k$ and the properties of our filtration function. 
The proof of assertion 4) is hence complete, because the diagram 
\[
 \xymatrix{
   H^{p_1}(E_k^{*,\overline{p_1}}) \otimes H^{p_2}(E_k^{*,\overline{p_2}}) \ar[r]^-{{\rm mult.}} 
   \ar[d]^{\beta \otimes \beta} & 
   H^{p_1+p_2}(E_k^{*,\overline{p_1} + \overline{p_2}}) \ar[d]^{\beta} \\
         H^{p_1}(X_G;\F_l) \otimes H^{p_2}(X_G;\F_l) \ar[r]^-{\cup}    &    H^{p_1 +p_2}(X_G;\F_l) 
 }
\]
commutes. 

Now assume that $H^*(X_G ; \F_l)$ vanishes in degrees $k$ and $k+1$. Let $\mu \in E^*$ be a
monomial in $t_1, \ldots, t_r$ of cohomological degree at least $k$. 
Because the elements $t_1, \ldots, t_r$ are cocycles in $E^2$ and 
because the differential acts as a derivation, it is enough to assume that $\mu$ has cohomological 
degree $k$ or $k+1$. If the cohomological degree is $k$, then $\mu$ represents the zero class in $H^*(E)$, 
because $H^k(E^*) \cong H^k(X_G ; \F_l) = 0$. If the cohomological degree is $k+1$, then we 
argue as follows. The elements $t_1, \ldots, t_r \in E^2$ correspond to cocycles 
$\overline{t_1} , \ldots, \overline{t_r} \in V_0^{2,1} \otimes \F_l$. 
By the choice of our filtration function, we obtain a 
monomial $\overline{\mu} \in E^{k+1, \overline{k}} \otimes \F_l$ corresponding to $\mu$. Because 
$\beta : H^{k+1}(E_k^{*,\overline{k}} \otimes \F_l) \to H^{k+1}(X_G ; \F_l)=0$ is injective, 
there is an element $\nu \in E_k^{k, \overline k} \otimes \F_l$ with $\overline{\mu} = d ( \nu)$, where $d$ is 
the differential on $E_k^{*,\overline k} \otimes \F_l$. 
But $\nu \in Z^* \otimes \F_l$ by construction of $Z^*$. Because $Z^* \otimes \F_l \to E^*$ is a cochain 
map, we conclude that $\mu$ is a cocycle in $E^*$, which finishes the proof of claim 5). 
\end{proof} 

If instead we start with a $G$-space $X$ which is not $l$-tame, then Theorem \ref{smallcochain} can 
still be applied in many cases. To illustrate this let $X$ be a connected nilpotent Kan complex 
with abelian fundamental group and let $G$ act on $X$. 
We denote by $X_{(l)}$ the $l$-localization of $X$, see \cite{BK}. By the 
assumptions on $X$ we have 
$\pi_t(X_{(l)}) = \pi_t(X) \otimes \Z_{(l)}$ for $t \geq 1$. Furthermore, 
by the functoriality of the $l$-localization, $X_{(l)}$ is again a $G$-space and the 
canonical $G$-equivariant map $X \to X_{(l)}$ induces an isomorphism 
\[
   H^*((X_{(l)})_G ; \F_l) \to H^*(X_G ; \F_l)  
\]
by a spectral sequence argument. Hence, if $X_{(l)}$ is $l$-tame, 
then the preceding discussion applies to $(X_{(l)})_G$ and for any $ k \geq 1$ with 
$\overline{k}< l$ we obtain a commutative graded
differential algebra $E$ over $\F_l$ which calculates the $\F_l$-cohomology of $X_G$ up to  degree $k$. 

We specialize this discussion to the $G$-space $X = S^{n_1} \times \ldots \times S^{n_k}$ appearing in Section \ref{cool}. Let  
$l > 3 \dim X$.  For each $ t \geq 1 $ with $\overline t < l$ this  assumption implies
\[
     \pi_t(S^{n_j}) \otimes \Z_{(l)} = \begin{cases} \Z_{(l)}  {\rm~for~} j > k_e {\rm~and~} t = n_j  \\
                                                     \Z_{(l)} {\rm~for~} j \leq k_e {\rm~and~} t = n_j {\rm~or~} t = 2n_j -1  \\
                                                      0 {\rm~otherwise} 
  \end{cases}
\]
for all $j \in \{1, \ldots, k\}$. Hence $\pi_t(X_{(l)})$ is a finitely generated free $\Z_{(l)}$ module of 
rank at most $k$. Because $l > 3 \dim X$, we have $l -1 > k$. This implies that  $G$ acts 
trivially on $\pi_t(X_{(l)})$, because the rational representation theory \cite{CR} of the group $\Z/l$ 
tells us that the smallest dimension of a nontrivial rational $\Z/l$-representation is $l-1$. 
We conclude that $X_{(l)}$ is an $l$-tame $G$-space in the sense of Definition \ref{tamedef}.   

Theorem \ref{smallcochain} and the subsequent discussion yield an $\F_l$-cochain algebra $E^*$ 
whose cohomology is isomorphic to $H^*(X_G ; \F_l)$ in degrees less than or equal to $\dim X +1$
and with the discribed behaviour of monomials in $t_1, \ldots, t_r$. This uses $l >  3 \dim X = \overline{\dim X + 1}$ 
so that we can  work with $k: = \dim X +1$. 

It follows from Theorem \ref{smallcochain} and the  partial calculation of 
$\pi_*((S^{n_j})_{(l)})$ displayed above 
that $E^*$ has the form described in  Proposition \ref{later}.

\section{Appendix: Proof of the tame Hirsch lemma} \label{append} 

The material in this appendix is based on Sections II.7 and II.8 of the diploma thesis of S\"orensen \cite{So}. 
Because the proof of the tame Hirsch lemma, a result  which is fundamental for our approach, 
is quite delicate, we provide the following detailed exposition.

Theorem \ref{tamehirsch} is first proven under the additional assumption 
that the unit 
\[
   \Psi_{\Lambda V}   : \Lambda V \to T(\|\Lambda V\|) 
\]
is a primary weak equivalence satisfying condition $k^+$. We will show in Proposition \ref{hui} below that there is no 
loss of generality in assuming this.

We work in several steps. At first, we consider the untwisted case, i.e. when $f$ 
is a constant map. Then  $E_f = X \times \| \Lambda V \|$ and $\Gamma_f$ is 
 induced by the projections $X \times \|\Lambda V\| \to X$, $X \times \|\Lambda V\| \to \| \Lambda V\|$, 
the unit $\Psi_{\Lambda V}$ and the product of forms in $T(X \times \| \Lambda V\|)$.

Now consider the commutative diagram of filtered graded differential algebras 
\begin{equation} \label{diagram} 
\xymatrix{T(X) \otimes \Lambda V \ar[d] \ar[r]^-{\id \otimes \Psi_{\Lambda V}} & T(X) \otimes T(\|\Lambda V \|)  \ar[d] \ar[r]^-{\wedge} & T(X \times \| \Lambda V \| ) \ar[dd]^{\rho}  \\
T(X) \botimes\Lambda V \ar[r]^-{\id  \botimes\Psi_{\Lambda V}} &    T(X) \botimes T(\| \Lambda V \|)\ar[d]^{\rho} \\
                                                      &   \big(  (T \botimes C)(X) \big)   \botimes \big( (T \botimes C)(\|\Lambda V\|) \big)  \ar[r]^-{\wedge \botimes \times}  
 & (T \botimes C)(X \times \| \Lambda V \|) \\
                                                      &   C(X) \botimes C(\| \Lambda V \|) \ar[u]_{\sigma}  \ar[r]^-{\times} & C(X \times \| \Lambda V \| ) \ar[u]_{\sigma} 
}
\end{equation}
where the maps $\rho$ and $\sigma$ are taken from the zig-zag sequence (\ref{zigzag}) and the bar indicates the filtrationwise 
tensor product.  

\begin{lem} \label{raffiniert} The canonical map 
\[
   T(X) \otimes\Lambda V \to T(X)  \botimes\Lambda V
\]
is a primary weak equivalence satisfying condition $k^+$. 
\end{lem} 

\begin{proof} For $n \geq 0$ let $(\Lambda^{(n)} V) \subset \Lambda V$ be the graded $\Q$-submodule 
generated by monomials consisting of exactly $n$ factors in $V$. Note that  $\Lambda^{(n)} V$  
inherits a filtration from $\Lambda V$ and that  $(\Lambda^{(n)} V)^{*,q} = 0$ for
$q <n \overline k$ and $(\Lambda^{(n)} V)^{*,q} = (\Lambda^{(n)} V)^{*,n \overline k} \otimes \Q_{q}$ for 
$q \geq  n \overline k$. Furthermore, the differential on $\Lambda V$ restricts 
to a differential on $\Lambda^{(n)} V$. 

For  $q \geq 1$ we have direct sum  decompositions
\begin{eqnarray*}  
(T(X) \otimes (\Lambda V))^{*,q} & = & \bigoplus_{n \geq 0 \, , ~ n \overline{k} \leq q}
   T^{*, q-n\overline{k}}(X) \otimes (\Lambda^{(n)} V)^{*, n \overline k} \otimes \Q_q  \\
(T(X) \botimes(\Lambda V))^{*,q} & = & \bigoplus_{n \geq 0 \, , ~n \overline{k} \leq q} T^{*,q}(X) \otimes(\Lambda^{(n)} V)^{*,q} \, . 
\end{eqnarray*}
It is therefore enough to investigate the canonical map 
\[
    T^{*,q-n \overline{k}}(X) \otimes (\Lambda^{(n)} V)^{*, n \overline k } \otimes \F_l \to T^{*,q}(X) \otimes
    (\Lambda^{(n)}V)^{*,q} \otimes \F_l
\]
for each $n \geq 0$ with $n \overline k \leq q$ and any prime $l > q$.  
By the algebraic K\"unneth formula for chain complexes over  the field $\F_l$, the induced 
map in cohomology splits in degree $t \geq 0$ into a direct sum of maps  
\[
\omega : H^i\big(T^{*,q-n \overline k }(X) \otimes\F_l\big) \otimes H^j\big((\Lambda^{(n)} V)^{*,n \overline k}
\otimes  \F_l \big)    \to  
H^i\big( T^{*, q}(X) \otimes \F_l\big) \otimes  H^j\big((\Lambda^{(n)} V)^{*, q} \otimes \F_l\big) \, ,
\]
where $i,j \geq 0$ and $i + j = t$. We have 
$H^j((\Lambda^{(n)} V)^{*, n \overline k} \otimes \F_l) = H^j((\Lambda^{(n)} V)^{*,q} \otimes \F_l)$,
 because $l >  q \geq n\overline k$. 
Now, if $q > n \overline{k}$, then  the map $H^*(T^{*,q-n\overline k}(X) \otimes \F_l) \to 
H^*(T^{*,q}(X) \otimes \F_l)$ is an isomorphism by the addendum to the Cenkl-Porter Theorem \ref{cenkl} 
and  hence $\omega$ is an isomorphism. 

We study the remaining case $q = n \overline{k}$.  Because $H^*((\Lambda^{(n)} V)^{*, n \overline k} \otimes \F_l)$ and 
$H^*((\Lambda^{(n)} V)^{*,q} \otimes \F_l)$ are concentrated in 
degrees at least $nd$, the map $\omega$ is an isomorphism if $i+j \leq nd$ 
and injective, if $i + j \leq nd+1$. Recall that for $i=1$ and $q = n \overline{k}$, 
the source of $\omega$ is the zero module.  

First let $t \geq 1$ and $q \geq \overline{t}$. This implies $n \geq 1$. 

If $n =1$, we get $t \leq nd$, because our assumption 
$\overline{t} \leq q = 1 \cdot \overline{k}$ implies $t \leq k$ and hence $t \leq d$ (recall that $d \geq k$ 
by assumption). 

If $n \geq  2$, then we have  $t+1 \leq nd$. For if we assume $t \geq nd$, then the obvious estimate 
$nd + ((t+1) - nd) \leq t + 1$  and the admissibility of 
the filtration function imply $n \overline{d} + \overline{(t+1)-nd} \leq \overline{t}$
and hence $n \overline{d} < \overline{t}$. If $d=1$, then this requires the assumption $n\geq 2$.
But this contradicts the assumed estimate $\overline{t} \leq q = n \overline{k} \leq n \overline{d}$.

These considerations show that $\omega$ is a $t$-equivalence, if $t \geq 1$ and $q \geq \overline{t}$. 

Now observe, that if $q \geq \overline{k}+1$, 
then the assumption $q = n \overline{k}$ automatically implies $n \geq 2$. Hence, the preceding 
estimate for the case $n \geq 2$ shows that 
under this assumption, $\omega$ is actually a $(t+1)$-equivalence. This completes the proof of 
Lemma \ref{raffiniert}. 
\end{proof} 

To continue the proof of Theorem \ref{tamehirsch} for constant $f$ we observe the following points.  Firstly, 
the map $\id \botimes\Psi_{\Lambda V}  : T(X) \botimes\Lambda V \to T(X)\botimes T(\|\Lambda V\|)$ is a primary weak equivalence 
satisfying condition  $k^+$. This follows from the K\"unneth formula \cite[Lemma 5.3.1]{Sp}
and the fact that by assumption $\Psi_{\Lambda V}$ is a primary weak equivalence satisfying condition  $k^+$. 
Secondly, the universal coefficient theorem \cite[Theorem 5.3.14]{Sp} implies that for each $q \geq 1$ and each prime $l > q$,
the maps $\rho$ and $\sigma$ induce isomorphisms in cohomology after tensoring the cochain complexes 
with $\F_l$. Thirdly, the cross product map 
\[ 
   C^*(X;\F_l) \otimes   C^*(\| \Lambda V \|;\F_l)  \to C^*(X \times \|\Lambda V \|; \F_l) 
\]
induces an isomorphism of cohomology groups for any prime $l$. This uses the K\"unneth formula for cohomology 
\cite[Theorem 5.6.1]{Sp} and the fact that for any prime $l$, 
the homology $H_*(\|\Lambda V\|;\F_l)$ is finitely generated in 
each degree. Recall that $\|\Lambda V \|$ is an Eilenberg-MacLane complex for a finitely generated 
$\Q_{\overline k}$-module, see Proposition \ref{eilen}.
This completes the proof of Theorem \ref{tamehirsch} for constant $f$.

The next step of the proof of Theorem \ref{tamehirsch} deals with the case when $f$ is homotopic to a constant map. But 
before we give a proof of this case, we will show the second main technical result of this section besides 
Lemma \ref{raffiniert}. 

\begin{lem} \label{aufwand} Let $A$ and $B$ be filtered CGDAs and  
let $V$ be an elementary complex of type 
$(d,\overline{k})$ with $d \geq k  \geq 1$. Assume that $\tau : V \to A$ is a filtered cochain map of 
degree $1$ and let $f : A \to B$ be a primary weak equivalence satisfying condition $k^+$ so that 
$f : A^{*,0} \to B^{*,0}$ is an isomorphism. Then the induced map 
\[
       A \otimes_{\tau}  \Lambda V \to B \otimes_{f \circ \tau} \Lambda V 
\]
is a primary weak equivalence satisfying condition $k^+$. 
\end{lem}

\begin{proof} For $q \geq 1$ let us consider the 
decreasing filtrations $\mathcal{F}$ of  $(A \otimes_{\tau} \Lambda V)^{*,q}$ and 
$(B \otimes_{f \circ \tau } \Lambda V)^{*,q}$ defined for $\gamma \geq 0$ by   
\begin{eqnarray*}
  \mathcal{F}_{\gamma} ( A \otimes_{\tau} \Lambda V)^{*,q} & = & (A^{\geq \gamma} \otimes_{\tau} \Lambda V)^{*,q} \\
  \mathcal{F}_{\gamma} ( B \otimes_{f \circ \tau} \Lambda V)^{*,q} & = & 
                          ( B^{\geq \gamma} \otimes_{f \circ \tau} \Lambda V)^{*,q} \, . 
\end{eqnarray*}
Let $l > q$ be a prime. The direct sum decomposition 
$\Lambda V = \bigoplus \Lambda^{(n)} V$ from the proof of Lemma \ref{raffiniert} leads to
the following direct sum decompositions of the $E_1$-terms of the resulting spectral sequences with $\F_l$-coefficients:  
\begin{eqnarray*}
   E_1^{i,j}(A) & = & \bigoplus_{n \geq 0\, , ~n \overline k \leq q} (A^{i,q - n \overline k} \otimes \F_l) 
        \otimes H^{j} ((\Lambda^{(n)} V)^{*,n \overline k} \otimes \F_l) \\
   E_1^{i,j}(B) & = & \bigoplus_{n \geq 0 \, ,~ n \overline k \leq q} (B^{i,q - n \overline k} \otimes \F_l) 
        \otimes H^{j} ((\Lambda^{(n)} V)^{*,n \overline k} \otimes \F_l)   \, . 
\end{eqnarray*}
On the $E_2$-terms we hence get the direct  sum decompositions 
\begin{eqnarray*}
E_2^{i,j}(A) & = & \bigoplus_{n \geq 0\, , ~ n\overline k \leq q} H^i (A^{*,q-n\overline k } \otimes \F_l) \otimes 
                     H^j (( \Lambda^{(n)} V)^{*,n \overline{k}} \otimes  \F_l)   \\
E_2^{i,j}(B) & = & \bigoplus_{n \geq 0\, , ~ n \overline k \leq q} H^i (B^{*,q-n\overline{k}} \otimes \F_l) \otimes 
  H^j ((  \Lambda^{(n)} V)^{*,n \overline{k}} \otimes \F_l)  \, , 
\end{eqnarray*}
but notice that the differentials $d_2$ on the $E_2$-terms do not in general
restrict to the single summands (with varying $i$ and $j$) unless $\tau$ is zero. 

We will now prove that the map of the $E_2$-terms induced by $f$ 
is a $t$-equivalence, if $t \geq 1$ and $q \geq \overline{t}$   and a $(t+1)$-equivalence, if 
$t \geq k+1$ and $q \geq \overline{t}$ or if $t = k$ and $q \geq \overline{k} +1$, where $t= i+j$ 
is the cohomological degree. This implies the statement of Proposition \ref{aufwand} by induction 
over the pages in the spectral sequence.   

First notice that $H^*( (\Lambda^{(n)} V)^{*, n \overline{k}} \otimes \F_l)$ is concentrated
in degrees at least $nd$ so that it suffices to show that 
\[
    \omega :  A^{*, q - n \overline{k}} \otimes  \F_l \to  B^{*, q-n \overline{k}} \otimes  \F_l
\]
is a $(t - nd)$-equivalence, respectively a $\big( (t+1)-nd\big)$-equivalence in the relevant cases. Since 
the map $\omega$ is itself an isomorphism for $q-n \overline{k}= 0 $ and the 
map $\omega$ is a $1$-equivalence for $q - n \overline{k} \geq 1 = \overline{1}$, the map $\omega$ is a $\big( (t+1) - nd\big) $-equivalence, 
if $-1 \leq \big( (t+1) - nd \big)  \leq 1$. Hence we will assume from now on that $t - nd \geq 1$. Because 
the case $n=0$ is trivial, we will also assume that $n\geq 1$. 
 
We start by analyzing the case $t \geq 1$ and $q \geq \overline{t}$. 

We will show at first that $\omega$ is a $(t - nd)$-equivalence. 
Because $t - nd \geq  1$ we can use the admissibility of the filtration function 
and the inequality $1 + nd + (t - nd) \leq t +1$  to conclude $\overline{1} + n \overline{d} + \overline{t-nd} \leq \overline{t}$ 
which implies $q - n \overline{k} \geq \overline{t - nd}$ 
(remember $k \leq d$). Hence $\omega$ is a $(t-nd)$-equivalence 
by the assumption on $f$. 

If in addition $n \geq 2$ or if $n=1$ and $d > 1$, then the inequality $nd + ((t+1)-nd) \leq t+1$ and the admissibility of 
the filtration function imply $q - n \overline{k} \geq \overline{(t+1)-nd}$  and $\omega$ is 
even a $\big( (t+1)-nd \big)$-equivalence.

Now let $n=1$ and $d=1$. This implies $k=1$. We will show that $\omega$ is still a $\big( (t+1)-nd \big)$-equivalence, 
if $t \geq k+1$, $q \geq \overline{t}$, 
or if $t = k$, $q \geq \overline{t}+1$. 

Let us first concentrate on the case $t \geq k+1$ and $q \geq \overline{t}$. 
If $t \geq k+2$, then 
$q - n\overline{k} = q -1 \geq \overline{t} -1 \geq \overline{t-1}$ and $t-1 \geq k+1$, so that 
$\omega$ is a $\big( (t-1)+1 \big)$-equivalence by the assumption on $f$. In other words, $\omega$ is 
a $\big( (t+1) - nd \big)$-equivalence as desired. 
If $t = k+1$, then (as $q \geq \overline t$) we get $q - n \overline{k} \geq \overline{k+1} - 1 = \overline{k} +1$ (remember $k=1$), so 
that $\omega$ is a $(k+1)$-equivalence by assumption on $f$ and hence $\omega$ is a $\big( (t+1)-nd \big)$-equivalence.

In the case $t = k$, $q \geq \overline{t} + 1$, we get $q - n \overline{k} = q-1 \geq \overline{t}$ so that 
$\omega$ is a $t$-equivalence by the assumption on $f$ and hence $\omega$ is again a $\big( (t+1) - nd \big)$-equivalence. 

This finishes the proof of Lemma \ref{aufwand}. 
\end{proof}

In order to proceed with the proof of Theorem \ref{tamehirsch}, we assume $f \simeq \const$ and 
choose a simplicial map  $H : X \times I \to \|\Lambda(\Sigma V)\|$ 
with $H \circ i_0 = f$ and $H \circ i_1 = \const$, where $I$ is the simplicial interval and 
$i_0, i_1: X \to X \times I$ are the canonical inclusions. 
For $\epsilon = 0$ and $\epsilon = 1$, we obtain a commutative diagram 
\[
\xymatrix{ 
    E_{\epsilon} \ar[r]^-{I_{\epsilon}} \ar[d] &  E_H  \ar[r]  \ar[d]  & \| \Lambda ( \cone V) \| \ar[d]\\
     X         \ar[r]^-{i_{\epsilon}}         &   X \times I \ar[r]^-{H}  & \| \Lambda( \Sigma V ) \| 
}  
\]
where the horizontal maps of the first square are {\em strong equivalences}. By definition 
this means that after applying the 
Cenkl-Porter functor $T$, one gets isomorphisms in cohomology in each degree and in each filtration level. 
Using the naturality of the construction of the map $\Gamma$, we hence get an induced commutative diagram 
\[
 \xymatrix{
   T(X) \otimes \Lambda(V) \ar[r]^-{\Gamma_{\const}}  &  T(E_{const}) \\
   T(X \times I) \otimes_{H^{\sharp}} \Lambda(V) \ar[r]^-{\Gamma_{H}} \ar[d]_{T(i_0) \otimes \id} \ar[u]^{T(i_1) \otimes \id} &  
   T(E_{H}) \ar[d]_{T(I_0)} \ar[u]^{T(I_1)} \\
   T(X) \otimes_{f^{\sharp}} \Lambda(V) \ar[r]^-{\Gamma_f} &       T(E_f) \\
  } 
\]
For any $q \geq 1$ and $l > q$ the vertical maps in this diagram restricted to filtration level $q$ 
induce isomorphisms of cohomology groups with $\F_l$-coefficients. 
This is clear for the right hand column. For the left hand column one argues with the spectral 
sequences introduced at the beginning of the proof 
of Lemma \ref{aufwand}. 

Hence, because $\Gamma_{\const}$ is a primary weak equivalence satisfying condition $k^+$, the same is true for the 
map $\Gamma_f$. 

Now these special cases are used to show the tame Hirsch lemma for arbitrary simplicial sets $X$
by a Mayer-Vietoris argument and induction over the $m$-skeleta of $X$. 

We consider simplicial maps  $f : X \to \| \Lambda( \Sigma V) \|$ where $X$ is a simplicial set and 
$V$ is an elementary complex of type $(d,\overline k)$ with $ d \geq k \geq 1$. In the following  
we leave $V$ fixed, but consider varying $X$. 

First, the tame Hirsch lemma holds, if $X$ is a disjoint union of simplicial $m$-simplices, $m \geq 0$. 
This follows from the preceding line of arguments, because in this case $f$ is homotopic to a constant 
map. Recall that $\| \Lambda (\Sigma V) \|$ is connected and fulfills the Kan condition by Proposition \ref{eilen}.
In particular, the tame Hirsch lemma holds, if $X$ is a disjoint union of $0$-dimensional simplicial sets. We now assume 
that the tame Hirsch lemma has been proven for any $(m-1)$-dimensional simplicial set $X$, where $m \geq 1$. If $X$ is $m$-dimensional, we write $X$ as a push out 
\[
      \xymatrix{  \coprod \partial \Delta[m] \ar[r]  \ar[d] & \coprod \Delta[m] \ar[d] \\
                  X^{m-1} \ar[r] &       X 
 }
\]
and do the same for $E_f$. Here $\partial \Delta[m]$ denotes the simplicial $(m-1)$-sphere. 
For  $q \geq 1$ we consider 
the following commutative diagram, where the columns are Mayer-Vietoris sequences and the horizontal arrows are induced 
by $\Gamma_f$ or appropriate restrictions thereof:
\begin{small}
\begin{equation} \label{MV} 
  \xymatrix{  0   \ar[d]       &                                                     0                        \ar[d]                         \\
   \left( T(X) \otimes_{f^{\sharp}} \Lambda V \right)^{*,q} \ar[d]   \ar[r] &   T^{*,q}(E_f )    \ar[d]      \\
   \left( T(X^{m-1}) \otimes_{f^{\sharp}} \Lambda V \right)^{*,q}  \oplus \left( T(\coprod \Delta[m]) \otimes_{f^{\sharp}} \Lambda V \right)^{*,q}  
       \ar[d] \ar[r]  &  T^{*,q}(E_{f}|_{X^{m-1}}) \oplus  T^{*,q}(E_{f}|_{\coprod  \Delta[m]})  \ar[d] \\
\left( T(\coprod \partial \Delta[m]) \otimes_{f^{\sharp}} \Lambda V \right)^{*,q} \ar[d]  \ar[r] &   T^{*,q}(E_{f}|_{\coprod \partial \Delta[m]})       \ar[d]   \\       0 & 0  } \, . 
\end{equation}
\end{small}
The right-hand column is exact because $T^{*,q}$ is a 
cohomology theory for $q \geq 1$. 
Unfortunately, this is not quite true for the left-hand column: indeed, for any simplicial set $Y$, 
we have a direct sum decomposition of graded $\Q_q$-modules
\begin{equation} \label{directsum}
    (T(Y) \otimes (\Lambda V))^{*,q}  =   \bigoplus_{n \geq 0 \, , ~ n \overline{k} \leq q}
    T^{*, q-n\overline{k}}(Y) \otimes (\Lambda^{(n)} V)^{*, n \overline k}   \otimes \Q_q
\end{equation}
as in the proof of Lemma \ref{raffiniert}. Furthermore,  if $q - n \overline{k}  \geq 1$  we have a short exact Mayer-Vietoris sequence
\[
     0  \to T^{*,q-n\overline k }(X)  \to T^{*,q-n\overline k}(X^{m-1})  \oplus  T^{*,q-n\overline k}\left( \coprod \Delta[m] \right)  \to T^{*,q-n\overline k}\left(\coprod \partial \Delta[m] \right)  \to 0 \, ,  
\]
which stays exact after tensoring with $(\Lambda^{(n)} V)^{*, n \overline k}$. However, for $q = n \overline k$ we
need to be careful: the sequence 
\[
     0  \to T^{*,0}(X)  \to T^{*,0}(X^{m-1})  \oplus  T^{*,0}\left( \coprod \Delta[m] \right)  \to T^{*,0}\left(\coprod \partial \Delta[m] \right)  
\]
is indeed exact for all $m \geq 1$, but  the final map is surjective, only if $m\geq 2$, in general. 
For any $q' \geq 0$ let us therefore consider the cochain complex  
\[
    \overline{T}^{*,q'}  \left(\coprod \partial \Delta[m] \right) := \im \Big( \big( T^{*,q'} (X^{m-1}) \oplus  T^{*,q'}( \coprod \Delta[m]) \big) \to T^{*,q'} (\coprod \partial \Delta[m]) \Big) \, , 
\]
which coincides with $T^{*,q'}( \coprod \partial \Delta[m])$ if $q' \geq 1$ or if $m \geq 2$.  For any prime $l > q$ we get an inclusion of cochain complexes
\[
    \left( \overline{T}\left( \coprod \partial \Delta[m]\right) \otimes_{f^{\sharp}} \Lambda V \right)^{*,q} \otimes \F_l \hookrightarrow
    \left( T\left(\coprod \partial \Delta[m]\right) \otimes_{f^{\sharp}} \Lambda V \right)^{*,q} \otimes \F_l \, , 
\]
which is a  $(t-1)$-equivalence for $t \geq 1$ and $q \geq \overline t$ and is a $t$-equivalence 
for $t \geq k$ and $q \geq \max\{\overline t, \overline k +1 \}$. This is implied by the direct sum decomposition (\ref{directsum}) 
and an argument similar to the one for the case $q = n\overline k$ 
in the proof of Lemma \ref{raffiniert}. 
Furthermore, the left 
hand column in Diagram (\ref{MV}) becomes  exact after replacing $T(\coprod \partial \Delta[m]) \otimes_{f^{\sharp}} \Lambda V$ 
by $\overline T (\coprod \partial \Delta[m]) \otimes_{f^{\sharp}} \Lambda V$. 

After these preliminaries a five lemma argument together with the induction hypothesis proves that 
$\Gamma_f: T(X) \otimes_{f^{\sharp}} \Lambda (V) \to T(E_f)$ 
is a primary weak equivalence satisfying condition $k^+$, thus finishing the induction step.

Finally, we assume that $X$ is not of finite dimension. In this case we observe that the restriction maps  $\big( T^{*,q}(X) \otimes_{f^{\sharp}} \Lambda V \big) \otimes \F_l\to\big(  T^{*,q}(X^m) \otimes_{f^{\sharp}} \Lambda V\big) \otimes \F_l $ and $T^{*,q}(E_f) \otimes \F_l \to T^{*,q}(E_{f}|_{X^m}) \otimes \F_l $ induce isomorphisms in cohomology up to degree $m-1$ for any prime $l > q$.  
       
As a  first application of the tame Hirsch lemma we study Eilenberg-MacLane complexes. The following result 
says that the assumption on the unit $\Psi_{\Lambda V}$ upon which the above proof of the tame Hirsch lemma was based always holds.  

\begin{prop} \label{hui} Let $V$ be an elementary complex of type $(d,\overline{k})$ with 
$d \geq k \geq 1$. Then the unit 
\[ 
   \Psi_{\Lambda V} : \Lambda V \to T(\| \Lambda V \| ) 
\]
is a primary weak equivalence satisfying condition $k^+$. 
\end{prop} 

The proof starts with the following lemma, which formulates one of 
the guiding principles of tame homotopy theory.

\begin{lem} \label{soul} Let $d,q \geq 1$ and let $l$ be a prime with $l > q$. Let 
$\Lambda(v)$ be the free (unfiltered) graded commutative $\F_l$-algebra generated by a 
variable $v$ of degree $d$.
Then the algebra map  
 \[
   \Lambda(v) \to H^*( K(\Q_q,d) ; \F_l) 
\]
induced by sending $v$ to the canonical generator of $H^d(K(\Q_q,d) ; \F_l)$
is an isomorphism in degree less than or equal to $d + 2l- 3$
and is an isomorphism in all degrees, if $d = 1$. 
\end{lem} 

\begin{proof} The result is clear, if $d=1$, because then both sides are concentrated in degree $1$. 
For an arbitrary prime $l$, a classical result of Cartan and Serre says that  
$H^*(K(\Z, d) ; \F_l)$ is a free commutative graded algebra over $\F_l$ in one generator of degree $d$ and 
further generators of degrees at least $d + 2(l-1)$ as the first reduced Steenrod power operation 
for the prime $l$ raises degree by $2(l-1)$. 

The inclusion $\Z \hookrightarrow \Q_q$ 
induces an isomorphism $H^*(K(\Q_q,d); \F_l) \to H^*(K(\Z,d) ; \F_l)$, because $l > q$.
This finishes the proof of Lemma \ref{soul}. 
\end{proof} 
 
\begin{prop} \label{single} Proposition \ref{hui} holds if $V$ is concentrated in degree $d$ and
$V^{d,\overline{k}} \cong \Q_{\overline k}$. 
\end{prop} 

\begin{proof} We start by noting that for each $q' \geq \overline{k}$ and for any 
prime $l > q'$, the unit $\Psi_{\Lambda V}$ induces an isomorphism 
\begin{equation} \label{isom} 
     H^{d}((\Lambda V)^{*,q'} \otimes \F_l) \to H^d(T^{*,q'}(\| \Lambda V \|) 
     \otimes \F_l)
\end{equation} 
by Proposition \ref{adjointbasic} and Theorem \ref{cenkl}.

We first treat the case $d=1$ separately. 
In this case we have $k=1$. Let $q \geq 1$ and $l > q$. Then $q \geq \overline{k}$ 
and  the cohomology $H^*(\Lambda(V)^{*,q} \otimes \F_l)$ is concentrated in degree $1$ and isomorphic to 
$\F_l$. The same holds for  
$H^*( T^{*,q}( \| \Lambda(V) \|) \otimes \F_l)$ by the Cenkl-Porter theorem. 
Hence, $\Psi_{\Lambda V}$ is a primary weak equivalence satisfying condition $k^+$.

From now on, we assume $d \geq 2$. Let $t \geq 1$ and $q \geq \overline{t}$. We distinguish 
the cases $t < k$ and $t \geq k$. 
In the first case we have $t < d$ (recall $k \leq d$) and 
we will show that $\Psi_{\Lambda V}: (\Lambda V)^{*,q} \otimes \F_l \to 
T^{*,q}(\|V \|) \otimes \F_l$ is a $t$-equivalence for $l > q$. That the 
induced map in cohomology is an isomorphism up to degree $t$ is trivial, because 
$t <d$. The injectivity in degree $t+1$ 
is again trivial, if $q < \overline{k}$, because $(\Lambda V)^{t+1,q} = 0$ in this case. 
For $q \geq \overline{k}$ this injectivity follows from the isomorphism (\ref{isom}). 

It remains to study the case $t \geq k$. 
Let $l > q$ be a prime. 
We will show that $\Psi_{\Lambda V} : (\Lambda V)^{*,q} \otimes \F_l \to 
T^{*,q}(\|V \|) \otimes \F_l$ is a $(t+1)$-equivalence. This will finish the proof of Proposition
\ref{single}. 

It follows from the Cenkl-Porter theorem and from Lemma \ref{soul}  that for all $\overline k \leq q' \leq q$ the 
graded groups $H^*(T^{*,q'} ( \| \Lambda(V) \|) \otimes \F_l)$ and  
$\Lambda(v) \otimes \F_l$ are canonically isomorphic up to degree $d + 2l -3$, 
where $v$ denotes a variable of degree $d$.
Furthermore, these isomorphisms are compatible with multiplications  
on filtration levels $q'_1$ and $q'_2$ with $\overline k \leq q_1', q_2' \leq q$ and $q'_1 + q'_2 \leq q$.  
Note that the  assumptions $t \geq k$ and $l > q \geq \overline{t}$ imply $l > \overline{k}$ 
so that Lemma \ref{soul} can be applied.  
As  $d \geq 2$ and $l \geq 2$, we obtain the inequality $d +2l-3 \geq l+1$. And because 
$l > q \geq \overline{t} \geq t$, the  
isomorphism $H^*(T^{*,q'} ( \| \Lambda(V) \|) \otimes \F_l) \cong  
\Lambda(v) \otimes \F_l$ holds up to degree $t+2$. For a proof of the  assertion 
that $\Psi_{\Lambda V}$ is a $(t+1)$-equivalence, it remains to 
show that for all $n \geq 1$ with $nd \leq t+1$ we get the inequality
$n \overline k \leq q$. For then we can apply induction on $n$ starting with the isomorphism 
(\ref{isom})  with $q' = \overline k$. Notice that in degree $t+2$, we need not care about 
the case  $n \overline k > q$, because then $(\Lambda V)^{t+2,q} \otimes \F_l = 0$ so that the map in question is 
certainly injective.  For $n = 1$ the required inequality holds, because $\overline{k} \leq \overline{t} \leq q$. 
For $n \geq 2$ it follows from the chain of inequalites 
$n  \overline{k} \leq n \overline{d} \leq \overline{t} \leq q$
where the second inequality uses the admissibility of the filtration function and the assumptions $d,n  \geq 2$.

This finishes the proof of Proposition \ref{single}.  
\end{proof}

We can now assemble the proof of Proposition \ref{hui}. 

First, let $V$ be concentrated in degree $d$. We proceed by 
induction on the rank of $V^{d, \overline{k}}$. If this rank is equal to $1$, we apply Proposition 
\ref{single}. Now let $V^{d, \overline{k}}$ be a free $\Q_{\overline{k}}$-module of rank $n+1$ with 
generators $v_1, \ldots, v_{n+1}$. Applying the tame Hirsch lemma and Proposition \ref{single}, the map 
\[
   T(\|\Lambda(v_1, \ldots, v_n) \|) \otimes \Lambda(v_{n+1}) \to T(\| \Lambda(v_1, \ldots, v_n) \| \times \| \Lambda(v_{n+1}) \|)
     = T(\| \Lambda V \|) 
\]
is a  primary weak equivalence satisfying condition $k^+$. Here $\Lambda(v_{n+1})$ is the free filtered CGDA generated by $v_{n+1}$. 
Hence we need to show that the 
map 
\[
   \Lambda(v_1, \ldots, v_n) \otimes \Lambda(v_{n+1}) \to T( \|\Lambda(v_1, \ldots, v_n)\| ) \otimes \Lambda(v_{n+1}) 
\]
is a primary weak equivalence satisfying condition   $k^+$. But this follows from Lemma \ref{aufwand} 
and the induction hypothesis. The assumption on filtration level $0$ holds, because 
$\| \Lambda ( v_1, \ldots, v_n) \|$ is connected. 

We are now able to prove the general case of Proposition \ref{hui}. Let $(V,\eta)$ be an elementary complex of 
type $(d, \overline{k})$ with $d \geq k \geq 1$. We can split $V$ into two elementary 
complexes $V_0$ of type $(d,\overline{k})$ 
and $V_1$ of type $(d+1, \overline{k})$ which are concentrated in degrees $d$ and $d+1$, respectively. 
The differential $\eta$ in $V$ is then regarded as a map of 
filtered cochain complexes  
$\eta : \Sigma V_0 \to V_1$. Applying the simplicial realization functor yields a simplicial map 
\[
     \| \Lambda ( \eta) \| : \|\Lambda V_1 \| \to \| \Lambda (\Sigma V_0) \| \, . 
\]
Because the unit $\Psi_{\Lambda V_0}$ is a primary weak equivalence satisfying condition $k^+$, which has 
already been shown, the tame Hirsch lemma applies 
and we get a primary weak equivalence satisfying condition $k^+$
\[
     \Gamma_{\|\Lambda (\eta) \|} :  T(\|\Lambda V_1\|) \otimes_{\Psi_{\Lambda V_1} \circ \eta} \Lambda V_0 \to 
     T( E_{\|\Lambda (\eta) \|}) \, . 
\]
By an argument similar to the one used in the proof 
of Proposition \ref{eilen}, we see that $E_{\|\Lambda (\eta) \| } = \| \Lambda V \|$ and  
by Lemma \ref{aufwand}, the map 
\[
    \Lambda V = \Lambda V_1  \otimes_{\eta} \Lambda V_0  \to T( \|\Lambda V_1  \|) \otimes_{\Psi_{\Lambda V_1 \circ \eta }} 
    \Lambda V_0 
\]
is a primary weak equivalence satisfying condition $k^+$. This finishes the proof of Proposition \ref{hui}.

\end{document}